\documentclass[11pt,reqno]{amsart}

\usepackage{amsmath, amsfonts, amssymb, amscd, enumerate, graphics, color, scalefnt}
\theoremstyle{plain}
\newtheorem{theo}{Theorem}[section]

\theoremstyle{definition}

\theoremstyle{plain}

\theoremstyle{definition}

\renewcommand{\=}{:=}

\newcommand{\beq}{\begin{equation}}
\newcommand{\eeq}{\end{equation}}
\renewcommand{\a}{\alpha}
\renewcommand{\b}{\beta}

\renewcommand{\d}{\delta}

\newcommand{\f}{\varphi}
\newcommand{\g}{\gamma}

\renewcommand{\k}{\kappa}
\renewcommand{\l}{\lambda}

\renewcommand{\o}{\omega}
\newcommand{\q}{\vartheta}

\renewcommand{\t}{\tau}

\renewcommand{\L}{\Lambda}
\renewcommand{\O}{\Omega}

\newcommand{\bC}{\mathbb{C}}

\newcommand{\bR}{\mathbb{R}}

\newcommand{\gc}{\mathfrak{c}}
\newcommand{\gd}{\mathfrak{d}}

\renewcommand{\gg}{\mathfrak{g}}
\newcommand{\gh}{\mathfrak{h}}

\newcommand{\gl}{\mathfrak{l}}
\newcommand{\gm}{\mathfrak{m}}

\newcommand{\so}{\mathfrak{so}}

\newcommand\SO{\mathrm{SO}}

\newcommand{\cB}{\mathcal{B}}

\newcommand{\cD}{\mathcal{D}}
\newcommand{\cE}{\mathcal{E}}
\newcommand{\cF}{\mathcal{F}}
\newcommand{\cG}{\mathcal{G}}
\newcommand{\cH}{\mathcal{H}}
\newcommand{\cI}{\mathcal{I}}

\newcommand{\cL}{\mathcal{L}}

\newcommand{\cS}{\mathcal{S}}
\newcommand{\cT}{\mathcal{T}}
\newcommand{\cU}{\mathcal{U}}
\newcommand{\cV}{\mathcal{V}}

\newcommand{\p}{\partial}

\renewcommand{\square}{\kern1pt\vbox
{\hrule height 0.6pt\hbox{\vrule width 0.6pt\hskip 3pt
\vbox{\vskip 6pt}\hskip 3pt\vrule width 0.6pt}\hrule height0.6pt}\kern1pt}

\DeclareMathOperator\Aut{Aut\;}

\DeclareMathOperator\Ad{Ad}
\DeclareMathOperator\ad{ad}

\renewcommand\Im{\operatorname{Im}}

\newcommand{\Hom}{{\operatorname{Hom}}}

\newcommand{\wt}{\widetilde}
\newcommand{\wh}{\widehat}

\newcommand{\be}{\begin{equation}}
\newcommand{\ee}{\end{equation}}

\def\<#1,#2>{\langle\,#1,\,#2\,\rangle}
\newcommand{\arr}{\begin{array}{rlll}}
\newcommand{\ea}{\end{array}}
\newcommand{\bea}{\begin{eqnarray}}
\newcommand{\eea}{\end{eqnarray}}
\newcommand{\bean}{\begin{eqnarray*}}
\newcommand{\eean}{\end{eqnarray*}}

\def\sideremark#1{\ifvmode\leavevmode\fi\vadjust{
\vbox to0pt{\hbox to 0pt{\hskip\hsize\hskip1em
\vbox{\hsize3cm\tiny\raggedright\pretolerance10000
\noindent #1\hfill}\hss}\vbox to8pt{\vfil}\vss}}}

\newcounter{ssig}
\newcounter{ttig}
\title[ Structure equations of   Levi degenerate  CR   hypersurfaces]
{Structure equations of Levi degenerate \\ CR hypersurfaces of uniform type}

\author[Costantino Medori and Andrea Spiro]{Costantino Medori and Andrea  Spiro}

\date{\today}


\thanks{{\it Acknowledgments}. This research was partially supported by the Project MIUR Real and Complex Manifolds  and by GNSAGA of INdAM}

\address{
Costantino Medori,
Dipartimento di Matematica e Informatica,
Universit\`a di Parma,
Parma, Italy.}
\email{costantino.medori@unipr.it}
\address{
Andrea Spiro, Scuola di Scienze e Tecnologie, Universit\`a di Camerino, 
Camerino, 
Italy.}
\email{andrea.spiro@unicam.it}

\begin{document}
\begin{abstract}  We explicitly  determine the structure equations of 5-dimensional Levi 2-nondegenerate CR hypersurfaces, using our
 recently constructed canonical Cartan connection for this class of CR manifolds. We also give 
 an outline of the basic properties of absolute parallelisms and Cartan connections, together with a brief discussion
 of the   absolute parallelisms   for such CR manifolds existing in the literature.
\end{abstract}
\maketitle
\null \vspace*{-.50in}

\section{Introduction}
\setcounter{equation}{0}
Let $M$ be  a  $(2n+1)$-dimensional CR hypersurface,  that is a manifold endowed with a pair $(\cD, J)$ formed by
\begin{itemize}
\item[a)] a    distribution $\cD \subset TM $ of  codimension $1$,
\item[b)] a smooth family $J$ of complex structures $J_x: \cD_x \longrightarrow \cD_x$, 
satisfying the   integrability condition,  i.e.   the complex distribution $\cD^{10} \subset T^\bC M$ of  the $(+i)$-eigenspaces of
the $J_x$   is involutive.  
\end{itemize}
We recall that the Levi form of $M$  is defined as follows.  For  $x \in M$, let  $\q$ be a $1$-form on a neighbourhood $\cU$ of $x$ with  $\ker \q_y = \cD_y$  at each point $y$ of $\cU$. The {\it Levi form at $x$}  is the symmetric bilinear map 
$$\cL_x: \cD_x \times \cD_x \to \bR,\qquad  \cL_x(v, w)\= d\q_x(v, Jw), $$
which is well known to be $J_x$-invariant and independent on the choice $\q$,  up to a scalar multiple.  
 If  the dimension of   $\ker \cL_x$  is constant over $M$, we call the CR hypersurface    {\it of uniform type}.
The case $\dim \ker \cL_x  =  0$ occurs if and only if the distribution $\cD$ is contact and 
in this case $(M, \cD, J)$ is called  {\it Levi-nondegenerate}.   If $\cD$ is of uniform type with $\dim \ker \cL_x  > 0$ at all points,  we call  it  {\it uniformly Levi-degenerate}.  
\par  
The simplest examples  of uniformly Levi-degenerate CR hypersurfaces  are given by the
cartesian products  $\overline M \times S$ of  a Levi-nondegenerate CR hypersurface $(\overline M, \overline \cD, \overline J)$ 
and an $m$-dimensional complex manifold $(S, J^S)$. The  natural CR structure of $\overline M \times S$  is the pair $(\cD, J)$ defined by 
$$\cD_x\= \overline \cD_{\bar x} + T_s S ,\qquad J_x \= \bar J_{\bar x} \times J^S_s\ \text{for all}\ x = (\bar x, s) \in \overline M \times S.$$
 If a  CR hypersurface  is locally   CR equivalent  with  a cartesian product of this kind 
 around any point, we say that {\it it admits  local CR straightenings}. \par
Under appropriate uniformity assumptions on   the  CR structure,  any uniformly Levi degenerate CR hypersurface  $(M, \cD, J)$
is equipped with  a  nested sequence  of complex distributions 
  \beq \label{freeman's sequence} \ldots\subset  \cF_2 \subset \cF_1 \subset \cF_0  \subset   \cF_{-1}\subset \cF_{-2} =  T^\bC M,
\eeq 
in which   $\cF_{-1} \= \cD^{10}$ and all  other subdistribution  $\cF_{i}$, $i \geq 0$, are 
 inductively defined in  a special way that  implies  that
$$[\cF_i, \cF_j] \subset \cF_{i+j}\qquad \text{for each}\ i, j \geq - 2$$ 
(here, we assume that   $\cF_{i + j} \= \ T^\bC M$   if  $i + j \leq -2$). 
 This nested sequence of distributions necessarily stabilises 
 after a finite number of steps and, by a result of Freeman (\cite{Fr}), it  has the  following crucial property:  {\it $(M, \cD, J)$  admits  local CR straightenings  if and only if  the first stabilising  distribution $ \cF_k$,  that is such that    $\cF_{k+\ell} = \cF_{k}$ for all $\ell\geq 0$,    
 is non trivial}.   \par
 The uniformly Levi degenerate CR hypersurfaces with trivial   stabilising  distribution $ \cF_k$ (hence, with no CR straightenings) are called {\it Levi $(k+1)$-nondegenerate}. This  notion  extends the concept of 
  Levi nondegeneracy,  since the Levi  $1$-nondegenerate hypersurfaces are precisely  the Levi nondegenerate  hypersurfaces  in the usual sense. 
  \par
The smallest  dimension for a CR hypersurface to be uniformly Levi degenerate  and with  no CR straightenings is  $5$. By dimension counting, any such $5$-dimensional CR hypersurfaces is    {\it 2-nondegenerate}. For conciseness,   we  call the CR manifolds of this kind  {\it girdled CR manifolds}.\par
The class of girdled CR manifolds and the associated  equivalence problem has been the main object of  investigation in  several  recent papers.   
In particular, in  \cite{MS}  we proved the existence of a canonical Cartan connection for any girdled CR manifold, obtaining in this way 
 a solution to  the equivalence problem  and a complete set of invariants for this class of CR manifolds. Independently and with preprints posted almost at the same time, 
 Isaev and Zaitsev presented in \cite{IZ} an alternative solution, hence another set of  invariants,  for  the same equivalence problem. Isaev and Zaitsev's  solution  is however  not corresponding to 
 a Cartan connection. Shortly after, a third solution  and another set of invariants  has been given  by  Pocchiola in \cite{Po}. 
 \par
Due to this, in several occasions,  the following question has been posed:  {\it  Is there a way to compare  one to the other such  solutions to the equivalence problem of girdled CR manifolds?}
\par
 Having this question in mind,  in this paper we  newly present  our solution   to the equivalence problem for girdled CR manifolds,  in  a way that allows
 an immediate comparison with the other existing  solutions. 
More precisely,  we first  provide 
 a quick review of the notions of equivalence problems, absolute parallelisms and    Cartan connections. 
 The intention of such overview  is twofold:  to fix  unambiguously the  meaning of all terms of our discussion
 and  to  clarify  the main reasons of interests for   solutions to equivalence problems coming from  
 canonical Cartan connections.  We then describe  in detail the canonical Cartan connections of girdled CR manifolds introduced in \cite{MS}, 
giving the  explicit expressions   of the corresponding structure equations and making manifest  all  
 curvature restrictions that characterise such connections.  
\par
\medskip
\section{Equivalence problems and Cartan connections}
\setcounter{equation}{0}
\label{sect2.2}
 Let  $\cG$ be  a class of geometric structures, that is of pairs $(M, \cS)$ formed by  a manifold $M$ with some  geometric datum $\cS$ of  fixed type (as,   for instance, a Riemannian metric $g$, a  distribution $\cD$, a CR structure $(\cD, J)$, etc.).  
Given two geometric structures $(M, \cS)$, $(M', \cS')$ in $\cG$,  the {\it local equivalences  around points  $x \in M$ and $x' \in M'$}   are the local diffeomorphisms  $f: \cU \to \cU'$ between neighbourhoods $\cU$, $\cU'$ of $x$, $x'$, transforming  $\cS|_{\cU}$  into  $\cS'|_{\cU'}$.  The   {\it equivalence problem for  the class $\cG$} it is the query for   an algorithm that   establishes when, given two points, there exists a local equivalence  around such two points. \par
A standard  approach   to  such  problem  consists in  looking for constructions that  give for  each   $(M, \cS)$ in  $\cG$  a unique triple $(P, (X_i), \wt{ (\cdot)})$  made of:
\begin{itemize}
\item[i)]  a  bundle $\pi: P \to M$ over the manifold $M$; 
\item[ii)] an absolute parallelism $(X_i)$ on $P$, i.e.
  an ordered $N$-tuples of  vector fields $(X_1, \ldots, X_N)$ that gives a frame at each tangent space $T_u P$; 
\item[iii)] an operator  $\wt{ (\cdot)}$  which maps each  local equivalence $f: \cU \to M' $ into  a bundle diffeomorphism $\wt f:  \cV \subset P \to P'$
  that projects onto $f$, 
\end{itemize}
 such   that the following holds: {\it a local diffeomorphism $F:  \cV \subset P \to P'$
 between the bundles  $P,P'$ of  two structures $(M, \cS)$, $(M', \cS')$ Êin $ \cG$ maps
the associated parallelisms $(X_i)$, $(X'_i)$  one into the other  if and only if  $F = \wt f$ for some local equivalence $f$}.
\par
Triples   $(P, (X_i), \wt{(\cdot)})$  with this property are  called {\it canonical absolute parallelisms for the class $\cG$} and any algorithm  that provides canonical absolute parallelisms 
 solves  the  equivalence problem for $\cG$  in the following sense. \par
Any absolute parallelism $(X_i)$  is  
uniquely determined by  the $N$-tuple of its dual  $1$-forms $(\o^1, \ldots, \o^N)$ and a   local diffeomorphism 
 transforms one absolute parallelism  into another if and only it it transforms the corresponding dual coframe fields  one  into the other.
We now observe that
the  differentials  $d\o^i$ admit unique expansions of the form $ d \o^i = \sum_{j < k} c^i_{jk} \o^j \wedge \o^k$.
These  
are  the so-called {\it structure equations of the  parallelism $(X_i)$} and the functions $c^i_{jk}$ are the associated  {\it first order invariants}.  Note that the  invariants  $c^i_{jk}$  can be explicitly  determined     from the  vector fields $X_i$ by recalling that
 \beq\label{alternative}  c^i_{jk} = d \o^i(X_j, X_k) = - \o^i([X_j, X_k]) .\eeq
  Their differentials    
 have  the form  
 $d c^i_{jk} = c^i_{jk|\ell_1} \o^{\ell_1}$ and 
  the  functions   $c^i_{jk|\ell_1}$ are called {\it second order invariants}. Their differentials   
  $d c^i_{jk|\ell_1} = c^i_{jk|\ell_1\ell_2} \o^{\ell_2}$  define the {\it third order invariants}  $c^i_{jk|\ell_1\ell_2}$  and so on. 
 By a fundamental result of Cartan and  Sternberg, if appropriate constant rank  conditions hold,  there is an  $m_o$ such that all
 invariants of order $r \leq m_o + 1$ give a map
 $
  F^{(m_o)} \= (c^i_{jk} ,\ c^i_{jk|\ell_1},\   c^i_{jk|\ell_1\ell_2},\  \ldots \  ,c^i_{jk|\ell_1\ldots \ell_{m_o}}): P \to \bR^{N_o}
$, 
which  completely characterises the  pair $(P, (X_i))$ up local equivalences  (\cite{St}, Thm. VII.4.1;   see also \cite{Kb, PS, Sp}).  
So,   any question on existence of  equivalences  between canonical   absolute
 parallelisms   of    $\cG$ is  in principle completely solvable by studying  the  invariants of the parallelisms
up to some finite order. This is the reason why  any algorithm that provides canonical 
absolute parallelisms for  $\cG$
is  considered  as a  solution to the equivalence problem for Êthis class.
 \par
\smallskip
Solutions of this type to the  equivalence problems are usually not unique. 
For instance,   the  so-called  {\it $G$-structures of finite type}  admit canonical  absolute parallelisms, 
determined via   a finite number steps,  each of them based on choices of certain normalising conditions  (\cite{St, Kb, Ta1, Mo, AS, PS}).  Different  choices lead   to non-equivalent canonical absolute parallelisms, hence to distinct solutions to the same equivalence problem.
Other examples 
are provided by  the celebrated  absolute parallelisms of  Chern and Moser  (\cite{CM}) and  
  of Tanaka  (\cite{Ta2, Ta3})   for   the  Levi-nondegenerate CR hypersurfaces, whose
  first order invariants are 
actually  constrained  by  non-equivalent sets of  linear equations. 
There exists also  three  distinct solutions to the equivalence problems for the  elliptic and hyperbolic CR manifolds   of codimension two, 
 which have been determined in  \cite{EIS, SSl, SS}.   
 \par
 \smallskip
Amongst all  canonical  absolute parallelisms that one  might   associate 
 with the  structures of a given class,  there are sometimes some special ones that  correspond to   Cartan connections.   As we will shortly see, 
  parallelisms   of this kind have  several very important additional   features. 
\par
We recall that 
a  {\it  Cartan connection  on  a manifold $M$,   modelled on a homogeneous space $G/H$,}  is 
a pair $(P, \varpi)$,  formed by 
 a principal $H$-bundle $\pi: P \to M$  and a $\gg$-valued 1-form 
 $\varpi:T P \to \gg = Lie(G)$  such that:
\begin{itemize}
\item[(a)] for each $y\in P$, the map  $\varpi_y: T_y P \to  \gg$ is a linear isomorphisms and
$(\varpi_y)^{-1}|_{\gh}: \gh  \to T^{V}_y P$
is the standard  isomorphism,  given by the  right action of $H$ on $P$,   between $\gh = Lie(H)$ and 
the  tangent space $T^{V}_y P$ of the fiber, 
\item[(b)] $R^*_h \varpi = \Ad_{h^{-1}} \circ \varpi$  for any $h\in H$. 
\end{itemize}
Given a class of geometric structures $\cG$, a correspondence between  the structures in $\cG$ and  Cartan connections on the underlying manifolds, 
is called  {\it canonical} if there is an associated  bijection  between the
local  equivalences $f: \cU   \to \cU' $ between manifolds $M, M'$ of $\cG$ and  the local diffeomorphisms $\wt f: P|_{\cU} \to P'|_{\cU'}$  between 
the bundles of the associated Cartan connections $(P, \varpi)$,  $(P', \varpi')$, that satisfy the condition
$\wt f^* \varpi' = \varpi$. \par
Note that if there is a  canonical Cartan connection $(P, \o)$ for any  manifold $M$ of  $\cG$,  each  basis $(E^o_i)$  for  $\gg = Lie(G)$ determines a canonical absolute 
parallelism $(P, (E_i), \wt{(\cdot)})$, formed by the   bundle $P$ and the absolute parallelism $(E_i)$  given by  the vector fields 
\beq \label{absoluteparall} E_i|_u := \varpi_u^{-1}(E^o_i),\qquad\ u \in P.\eeq
Hence,  {\it  any construction of canonical  Cartan connections for a class $\cG$ automatically provides   a  solution to the corresponding  equivalence problem}. \par
However,  the interest for Cartan connections is by far much wider than  their uses for  equivalence problems.  For an introduction to  the variety  of  possible applications,  see  e.g.  \cite{Kb, Sh, CSS1, CSS2, CSS3, CG, SS1, BES, Wi} and  references therein. \par
One of the  most  basic reasons of interest for  Cartan connections is given by the following  fact:  {\it  If  $(P, \varpi)$  is a Cartan connection  on  $M$ modelled on $G/H$,  the associated  $\gg$-valued {\it curvature $2$-form}  $\O =  d \varpi  + \frac{1}{2} [\varpi, \varpi]$ on $P$
vanishes identically if and only if  $P$ is locally equivalent to the Lie group $G$ and $M$ is locally equivalent to the homogeneous model $G/H$.}
This means that if the elements of a  class $\cG$ of geometric structures admit canonical Cartan connections modelled on a given homogeneous spaces, for each of them  there exists a very informative   indicator (namely, the curvature $\O$) 
of  how it locally deviates from   the homogeneous model.  \par
From this and other facts on  Cartan connections,  one has also that  geometric structures  admitting canonical Cartan connections are
 equipped with   distinguished  families of  appropriate curves or submanifolds of higher dimension, which are invariant under  local equivalences and 
  play the same role of   geodesics  and  chains  in Riemannian geometry and in  geometry of   Levi non-degenerate hypersurfaces, respectively (see e.g. \cite{BDS, SSl}).
Such distinguished  curves and submanifolds can be also combined and determine systems of normal coordinates, which  allow to reduce several  questions  to  geometric properties  of the  homogeneous models (see, for instance,  \cite{SSl, SS1}). 
 \par
At the best of our knowledge, the first  methodical study on  the possibilities of  constructing canonical Cartan connections was done by Tanaka in \cite{Ta3}. There he proved the existence of canonical Cartan connections for  an important family of classes of geometric structures,  modelled on homogeneous spaces 
$G/H$ of (semi)simple Lie groups and  with parabolic isotropy subgroups  $H \subset G$.  His results
were later extended in various senses   by T. Morimoto in \cite{Mo} and \v Cap and Schichl in \cite{CS}. For a concise review of  Tanaka's  results, see \cite{AS}. \par
We conclude this short discussion of Cartan connections recalling that in \cite{AS}, Alekseevsky and the second author proved that Tanaka's method of construction of  Cartan connections can be considered as a derivation of  a more general method of construction of  absolute parallelisms, also invented by Tanaka (\cite{Ta1}).  
This second method applies to a wider range of geometric structures, called {\it Tanaka's structures of finite type}, and 
produces canonical parallelisms $(P, (X_i), \wt{(\cdot)})$,   formed by bundles  $\pi: P\to M$ that   {\it in general  are not principal bundles} and by  parallelisms
$(X_i)$ that  {\it in general  are not   determined by  a $\gg$-valued   $1$-form} $\varpi$  satisfying the properties of Cartan connections. Nonetheless, 
 for a special class of Tanaka structures,  modelled on  homogeneous spaces 
$G/H$ of a semisimple $G$ and  parabolic subgroup $H \subset G$, 
the general  construction  can be performed in such a way that it produces
 a bundle $\pi: P\to M$, which {\it is a  principal $H$-bundle}, and an absolute parallelism $(X_i)$ on $P$, 
  which {\it is determined  by  a  Cartan connection $\varpi$} (see \cite{AS} for details).   
\par
\medskip
\section{Cartan connections of girdled CR manifolds\\
and corresponding  structure equations} 
\setcounter{equation}{0}
Now we focus on girdled CR manifolds $(M, \cD, J)$, i.e. on  $5$-dimensional CR  hypersurfaces of uniform type, which are Levi $2$-non\-degenerate. 
As we already  mentioned in the introduction, any such CR hypersurface is  Levi degenerate and yet admits no local straightenings. The name {\it girdled}  has been chosen to  allude  to  such lack of  straightenings.\par
 \medskip
One of the most important  examples of  girdled CR manifolds and, as we will shortly see,  a model  for these geometric structures is given by the following homogeneous manifold. 
Consider
the bilinear form $(\cdot , \cdot )$ and the pseudo-Hermitian form  $<\cdot , \cdot>$ on $\bC^5$ defined by 
\beq \label{inner} (t, s) = t^T I_{3,2}  s,\qquad <t, s> = (\overline t, s) ,\qquad  I_{3,2} =  \left(\begin{array}{c|c} I_3 & 0\\ \hline 
0 & - I_2 \end{array} \right),\eeq
and the corresponding semi-algebraic subset   $M_o \subset \bC P^4$  defined by 
\beq\label{tubeflc1} \left\{\begin{array}{l}  (t, t) =  (t^0)^2 + (t^1)^2 + (t^2)^2 - (t^3)^2   -  (t^4)^2= 0,\\[4pt]
< t, t> =  |t^0|^2 + |t^1|^2 + |t^2|^2 - |t^3|^2   - |t^4|^2 = 0,\\[4pt]
\Im\left(t^3 \overline{ t^4}\right) >  0.
\end{array}\right.\eeq
 It is known (see e.g. \cite{SV}) that  $M_o$ 
  is a  $\SO_{3,2}^o$-homogeneous,  5-dimensional CR submanifold  of $\bC P^4$  (here, $\SO_{3,2}^o$ is the  identity component of $\SO_{3,2}$)   and contains  $\cT_o = M_o \cap \{ \Im(t^3 \overline{(t^0 + t^4)}) > 0\}$  as   open   dense  subset, which  is  CR equivalent to the
so called {\it tube over the future light cone in $\bC^3$}, i.e. to  the real hypersurface
\beq\label{tubeflc2} \cT = \{\,(z^1, z^2, z^3) \in \bC^3\,: \,(x^1)^2 + (x^2)^2 - (x^3)^2 = 0,\,x^3 > 0\,\}. \eeq
It  turns out  that  $M_o$  is   girdled and its group of CR automorphisms coincides with  $\Aut(M_o) = \SO_{3,2}^o$. 
Hence, if we 
denote by $H \subset \SO_{3,2}^o$ the isotropy subgroup of $\Aut(M_o) $ at  some point, 
$M_o$ is CR equivalent to  the homogeneous space $\SO_{3,2}^o/H$, equipped with an appropriate invariant girdled CR structure. \par
The homogeneous CR manifold $M_o =\SO_{3,2}^o/H$ is a  modelling space, of which  any girdled CR manifold can be considered as a local deformation.  This is a consequence of the main theorem of our paper  \cite{MS},  namely  
\begin{theo} \label{main} For any  5-dimensional girdled CR manifold  $(M, \cD, J)$,   there exists a canonical 
Cartan connection $(Q, \varpi)$,  modelled on  the homogeneous  CR manifold  $M_o = \SO_{3,2}^o/H$ described above.  
\end{theo}
The proof  of this theorem is constructive and provides an explicit description of the  bundle $\pi: Q\longrightarrow M$ and of the $\so_{3,2}$-valued 1-form $\varpi$ (more precisely, 
of a collection of vector fields, by which   $\varpi$  is uniquely determined). 
Our construction  is based on a modification of    Tanaka's general scheme  for building up   absolute parallelisms. The fact that our collection of  vector fields   actually defines 
 a Cartan connection is a consequence of an appropriate   tuning of  each step of the construction. 

 As    it is shown  in \cite{AS} (see  also above,  end of  \S \ref{sect2.2}), 
  even   the classical   Tanaka's  method   can be used to produce Cartan connections,  provided  that    appropriate algebraic conditions are satisfied. Such conditions  certainly occur  for the 
{\it parabolic geometries} \cite{CS}, i.e. 
the geometric structures   modelled on homogeneous spaces $G/H$ of   semisimple  Lie groups $G$  with parabolic   $H$. 
Since the girdled CR manifolds are modelled on a homogeneous space $G/H$ of the semisimple Lie group  $G = \SO^o_{3,2}$ 
 with  a {\it non parabolic}  $H$, 
our result shows  that the   above conditions 
might  occur  for a  wider  and interesting class of homogeneous models.\par
 \smallskip
As  pointed out in  \S \ref{sect2.2}, the absolute parallelism,  that is determined by  the canonical   Cartan connection $(Q, \varpi)$  and a basis of $\so_{3,2}$,   provides a  solution to the equivalence problem for  girdled CR manifolds.  At the best of our knowledge, at the moment there are  two  other absolute parallelisms for girdled CR manifolds, hence two other  solutions to the same problem (\cite{IZ, Po}), but  none of them corresponds to a  Cartan connection. 
\par
\smallskip
In the next sections, we select a special basis for $\so_{3,2}$ and we write explicitly the structure equations of the absolute parallelism corresponding to such special basis. 
Such explicit expressions also allow  immediate comparisons with the structure equations of  the  parallelisms provided by the other solutions to the   equivalence problem for girdled CR manifolds. \par
\smallskip
\subsection{A convenient  basis for $\so_{3,2}$}
The Lie algebra $\gg=\so_{3,2}$  has a natural structure of graded Lie algebra, which can be explicitly  described as follows. 
Consider  a   system of projective coordinates on $\bC P^4$, in which   the  scalar product $(\cdot, \cdot)$ defined in \eqref{inner} assumes the form
\beq \label{projcoor}(t, s) = t^T\cI\,s\qquad  \text{with}\qquad \cI = \hbox{$\scalefont{0.75}{\left(\arraycolsep=2pt
\begin{array}{cc|c|cc}0&0&0&0&1 \\[-0.5pt]
0&0&0&1&0\\[0.5pt]
\hline
0&0&1&0&0\\
\hline
0&1&0&0&0\\[-0.5pt]
1&0&0&0&0
\end{array}
\right)}$}.\eeq
By means of these new coordinates, the Lie algebra   $ \so_{3,2}$  of the isometries of $(\cdot, \cdot)$ can be identified with the  Lie algebra of  real matrices $A$ such that  $A^T \cI + \cI A = 0$,  i.e.,  of the form 
$$A = \hbox{$\scalefont{0.75}{\left(\arraycolsep=2pt
\begin{array}{c|c|c}
\begin{matrix} a_1 & a_2 \\a_3 & a_4 \end{matrix} &  \begin{matrix}  a_5 \\a_6 \end{matrix} & \begin{matrix}  a_7 & 0 \\0 & -a_7 \end{matrix}\\
\hline 
\begin{matrix}   a_{8} &  a_{9} \end{matrix} & 0  & \begin{matrix}  -a_{6} & -a_5 \end{matrix}\\
 \hline
\begin{matrix} a_{10} & 0 \\0 & -a_{10} \end{matrix} &  \begin{matrix}  - a_{9}\\-a_{8} \end{matrix} & \begin{matrix}- a_4  & - a_2 \\- a_3 & - a_1 \end{matrix}
\end{array}
\right)}$},\qquad\text{for some}\, \ a_i\in \bR.$$
This shows that $\so_{3,2}$   is  the direct sum of the vector subspaces 
\beq \label{gradedspaces}
\begin{split} &\gg_{-2} = < e^o_{-2} >,\,\gg_{-1} = < e^o_{-1|1}, e^o_{-1|2}>,\,  \gg_0 = < e^o_{0|1}, e^o_{0|2}, E^o_{0|1}, E^o_{0|2}> ,\\
&\gg_1 = < E^o_{1|1}, E^o_{1|2}>,\,\gg_2 = < E^o_2>,
\end{split}
\eeq
spanned by the matrices
\beq 
\begin{split}
\label{basisso32}
& 
e^o_{-2} {=}   \hbox{$\scalefont{0.75}{\left(\arraycolsep=2pt
\begin{array}{c|c|c}
\arraycolsep=2pt\begin{matrix} 0 & 0 \\[-3pt] 0 &0 \end{matrix} & \begin{matrix} 0 \\[-3pt] 0  \end{matrix} &  \begin{matrix}  0 & 0 \\[-3pt]0 & 0 \end{matrix}\\[-1.5pt]
\hline 
\begin{matrix}  0 &   0 \end{matrix} & 0  & \begin{matrix}  0 &   0 \end{matrix}\\[-1.5pt]
 \hline
\begin{matrix}  1 & 0 \\[-3pt]0 & -1 \end{matrix} &  \begin{matrix} 0 \\[-3pt] 0  \end{matrix} & \begin{matrix} 0  & 0 \\[-3pt]0 & 0 \end{matrix}
\end{array}
\right)} $}, \ e^o_{-1|1} {=}   \hbox{$\scalefont{0.75}{\left(\arraycolsep=2pt
\begin{array}{c|c|c}
\begin{matrix} 0 & 0 \\[-3pt]0 &0 \end{matrix} & \begin{matrix} 0 \\[-3pt] 0  \end{matrix} &  \begin{matrix}  0 & 0 \\[-3pt]0 & 0 \end{matrix}\\[-1.5pt]
\hline 
\begin{matrix}  1 &   0 \end{matrix} & 0  & \begin{matrix}  0 &   0 \end{matrix}\\[-1.5pt]
 \hline
\begin{matrix}  0 & 0 \\[-3pt]0 & 0 \end{matrix} &  \begin{matrix} 0 \\[-3pt] -1  \end{matrix} & \begin{matrix}0  & 0 \\[-3pt]0 & 0 \end{matrix}
\end{array}
\right)}$}, \ e^o_{-1|2} {=}  \hbox{$\scalefont{0.75}{\left(\arraycolsep=2pt
\begin{array}{c|c|c}
\begin{matrix} 0 & 0 \\[-3pt]0 &0 \end{matrix} & \begin{matrix} 0 \\[-3pt] 0  \end{matrix} &  \begin{matrix}  0 & 0 \\[-3pt]0 & 0 \end{matrix}\\[-1.5pt]
\hline 
\begin{matrix}  0 &   1 \end{matrix} & 0  & \begin{matrix}  0 &   0 \end{matrix}\\[-1.5pt]
 \hline
\begin{matrix}  0 & 0 \\[-3pt]0 & 0 \end{matrix} &  \begin{matrix} -1 \\[-3pt] 0  \end{matrix} & \begin{matrix}0  & 0 \\[-3pt]0 & 0 \end{matrix}
\end{array}
\right)}$},
\\
& 
e^o_{0|1} {=}   \hbox{$\scalefont{0.75}{\left(\arraycolsep=2pt\
\begin{array}{c|c|c}
\begin{matrix} 1 & 0 \\[-3pt]0 & - 1 \end{matrix} & \begin{matrix} 0 \\[-3pt] 0  \end{matrix} &  \begin{matrix}  0 & 0 \\[-3pt]0 & 0 \end{matrix}\\[-1.5pt]
\hline 
\begin{matrix}  0 &   0 \end{matrix} & 0  & \begin{matrix}  0 &   0 \end{matrix}\\[-1.5pt]
 \hline
\begin{matrix}  0 & 0 \\[-3pt]0 & 0 \end{matrix} &  \begin{matrix} 0 \\[-3pt] 0  \end{matrix} & \begin{matrix} 1  & 0 \\[-3pt]0 & -1 \end{matrix}
\end{array}\right)}$},
\ 
e^o_{0|2} {=}  \hbox{$\scalefont{0.75}{\left(\arraycolsep=2pt
\begin{array}{c|c|c}
\begin{matrix} 0 & 1 \\[-3pt]1 &0 \end{matrix} & \begin{matrix} 0 \\[-3pt] 0  \end{matrix} &  \begin{matrix}  0 & 0 \\[-3pt]0 & 0 \end{matrix}\\[-1.5pt]
\hline 
\begin{matrix}  0 &   0 \end{matrix} & 0  & \begin{matrix}  0 &   0 \end{matrix}\\[-1.5pt]
 \hline
\begin{matrix}  0 & 0 \\[-3pt]0 & 0 \end{matrix} &  \begin{matrix} 0 \\[-3pt] 0  \end{matrix} & \begin{matrix}0  & -1 \\[-3pt]-1 & 0 \end{matrix}
\end{array}
\right)}$},
\\
& 
E^o_{0|1} {=}   \hbox{$\scalefont{0.75}{\left(\arraycolsep=2pt
\begin{array}{c|c|c}
\begin{matrix} 1 & 0 \\[-3pt]0 &1 \end{matrix} & \begin{matrix} 0 \\[-3pt] 0  \end{matrix} &  \begin{matrix}  0 & 0 \\[-3pt]0 & 0 \end{matrix}\\[-1.5pt]
\hline 
\begin{matrix}  0 &   0 \end{matrix} & 0  & \begin{matrix}  0 &   0 \end{matrix}\\[-1.5pt]
 \hline
\begin{matrix}  0 & 0 \\[-3pt]0 & 0 \end{matrix} &  \begin{matrix} 0 \\[-3pt] 0  \end{matrix} & \begin{matrix}- 1  & 0 \\[-3pt]0 & -1 \end{matrix}
\end{array}
\right)}$},
\ 
E^o_{0|2} {=}   \hbox{$\scalefont{0.75}{\left(\arraycolsep=2pt
\begin{array}{c|c|c}
\begin{matrix} 0 & 1 \\[-3pt]-1 &0 \end{matrix} & \begin{matrix} 0 \\[-3pt] 0  \end{matrix} &  \begin{matrix}  0 & 0 \\[-3pt]0 & 0 \end{matrix}\\[-1.5pt]
\hline 
\begin{matrix}  0 &   0 \end{matrix} & 0  & \begin{matrix}  0 &   0 \end{matrix}\\[-1.5pt]
 \hline
\begin{matrix}  0 & 0 \\[-3pt]0 & 0 \end{matrix} &  \begin{matrix} 0 \\[-3pt] 0  \end{matrix} & \begin{matrix}0  & -1 \\[-3pt]1 & 0 \end{matrix}
\end{array}
\right)}$},\\
& E^o_{1|1}{=}   \hbox{$\scalefont{0.75}{\left(\arraycolsep=2pt
\begin{array}{c|c|c}
\begin{matrix} 0 & 0 \\[-3pt]0 &0 \end{matrix} & \begin{matrix} 1 \\[-3pt] 0  \end{matrix} &  \begin{matrix}  0 & 0 \\[-3pt]0 & 0 \end{matrix}\\[-1.5pt]
\hline 
\begin{matrix}  0 &   0 \end{matrix} & 0  & \begin{matrix}  0 &   -1 \end{matrix}\\[-1.5pt]
 \hline
\begin{matrix}  0 & 0 \\[-3pt]0 & 0 \end{matrix} &  \begin{matrix} 0 \\[-3pt] 0  \end{matrix} & \begin{matrix}0  & 0 \\[-3pt]0 & 0 \end{matrix}
\end{array}
\right)}$},\ 
E^o_{1|2} {=}  \hbox{$\scalefont{0.75}{\left(\arraycolsep=2pt
\begin{array}{c|c|c}
\begin{matrix} 0 & 0 \\[-3pt]0 &0 \end{matrix} & \begin{matrix} 0 \\[-3pt] 1  \end{matrix} &  \begin{matrix}  0 & 0 \\[-3pt]0 & 0 \end{matrix}\\[-1.5pt]
\hline 
\begin{matrix}  0 &   0 \end{matrix} & 0  & \begin{matrix}  -1 &   0 \end{matrix}\\[-1.5pt]
 \hline
\begin{matrix}  0 & 0 \\[-3pt]0 & 0 \end{matrix} &  \begin{matrix} 0 \\[-3pt] 0  \end{matrix} & \begin{matrix}0  & 0 \\[-3pt]0 & 0 \end{matrix}
\end{array}
\right)}$},\ 
E^o_2 {=}   \hbox{$\scalefont{0.75}{\left(\arraycolsep=2pt
\begin{array}{c|c|c}
\begin{matrix} 0 & 0 \\[-3pt]0 &0 \end{matrix} & \begin{matrix} 0 \\[-3pt] 0  \end{matrix} &  \begin{matrix}  1 & 0 \\[-3pt]0 & -1 \end{matrix}\\[-1.5pt]
\hline 
\begin{matrix}  0 &   0 \end{matrix} & 0  & \begin{matrix}  0 &   0 \end{matrix}\\[-1.5pt]
 \hline
\begin{matrix}  0 & 0 \\[-3pt]0 & 0 \end{matrix} &  \begin{matrix} 0 \\[-3pt] 0  \end{matrix} & \begin{matrix} 0  & 0 \\[-3pt]0 & 0 \end{matrix}
\end{array}
\right)}$}.
\end{split}
\eeq 
We will refer to  the collection $\cB^o$ of these matrices  as  {\it standard basis} of $\so_{3,2}$. \par
 Note  that 
the  each $\gg_k$ in \eqref{gradedspaces} is  the  eigenspace
 of the adjoint action of the {\it grading element} $Z\= E^o_{0|1}$  with  eigenvalue $k$, so that 
 $$[\gg_i, \gg_j] \subset \gg_{i+j}\qquad \text{for all}\ i, j,$$
 where, by convention,  we assume  $\gg_k = \{0\}$ for any  $k  \notin \{-2,-1, 0, 2, 2\}$. In other words,  {\it $\so_{3,2}$ has a natural structure of graded Lie algebra}.  
\par
We note  that also the  Lie algebra $\gh = Lie(H)$ of the isotropy subgroup $ H \subset \SO^o_{3,2}$ at  $x_o = [1 : i : 0: 0: 0]$ is natural graded. Indeed,  it  decomposes into the direct sum
$\gh = \gh_0 + \gg_1 + \gg_2$ with $\gh_0\= \langle E_{0|1}, E_{0|2} \rangle$.
\par
\smallskip
We conclude this section introducing another convenient basis for $\so_{3,2}$, which is a technical  modification of $\cB^o$,  more suitable 
for several  arguments concerning the CR structure of  the model space $\SO^o_{3,2} {\cdot} x_o \simeq  \SO^o_{3,2}/H$. 
Indeed, in many places it is   more appropriate to  consider instead of  the  four elements $E^o_{\ell|j}$, $e^o_{-\ell|j}$,    $\ell = 0,1$,  $j = 1,2$,     the four complex matrices  in  $(\so_{3,2})^\bC$:
\beq \label{definitionholelements}
\begin{array}{ll}
 E^o_{\ell (10)} =  \frac{1}{2} \left(E^o_{\ell|1} - i E_{\ell|2}\right), &
E^o_{\ell(01)} = \overline{E^o_{\ell(10)}}
,\\[6pt]
e^o_{-\ell(10)} =  \frac{1}{2} \left(e^o_{-\ell|1} - i e^o_{-\ell|2}\right)
 , &e^o_{-\ell(01)} = \overline{e^o_{-\ell(10)}},
 \end{array} \quad \ell = 0,1.\eeq
So, in the following, instead of expanding 
 the  elements of $\so_{3,2}$ in terms of the standard basis $\cB^o$,  we often expand the same elements 
in terms of the {\it  standard CR basis} 
\beq \label{modifiedbasis} \cB^{CR} = \left(e^o_{-2}, e^o_{-1(10)}, e^o_{-1(01)}, e^o_{0 (10)}, e^o_{0 (01)}, E^o_{0 (10)}, E^o_{0 (01)}, E^o_{1 (10)}, E^o_{1(01)}, E^o_{2}\right).\eeq
Since  the elements  $ X \in \so_{3,2}$ are    real matrices, 
 their  expansion  in the    standard  CR basis have the form 
 $X = \sum \l^A e^o_A  +  \sum \mu^A E^o_A$, with 
  coefficients satisfying the reality conditions
 $$\l^{-2}, \mu^{-2} \in \bR\qquad\text{and}\qquad \l^{-\ell(01)} = \overline{\l^{-\ell(10)} },\  \mu^{\ell(01)} = \overline{\mu^{\ell(10)} }\ \text{for }\ \ell =  0, 1.$$
A table of  all  Lie brackets between  elements   in   $\cB^{CR}$ can be found in \cite{MS}. \par
\medskip
\subsection{The absolute parallelism associated  with  the  standard basis}
\label{sect3.3}
Consider now a girdled CR manifold $(M, \cD, J)$ and its canonical Cartan connection $(Q, \varpi)$ modelled
on $M_o = \SO^o_{3,2}/H$. As we discussed in \S \ref{sect2.2}, the relation  \eqref{absoluteparall}  associates with  each element $e^o_A$ or $E^o_B$  of 
  the standard basis $\cB^o$ of $\so_{3,2}$ a vector field that we denote by $e_A$ or $E_B$, respectively. 
The ordered $10$-tuple $(e_A, E_B)$  is the {\it absolute parallelism  corresponding   to the  basis $\cB^o$}. \par
As we observed,  in place of these (real) vector fields,  it is often more convenient to consider  the 
 collection of (real and complex) vector fields 
\beq \label{complexabsoluteparall} \left(e_{-2}, e_{-1(10)}, e_{-1(01)}, e_{0 (10)}, e_{0(01)}, E_{0 (10)}, E_{0(01)}, E_{1(10)}, E_{1(01)}, E_{2}\right),\eeq
with   $e_{-\ell (10)}, e_{-\ell (01)}, E_{\ell (10)}, E_{\ell(01)}$, $\ell = 0,1$,     defined by 
\beq \label{definitionholelements-bis}
\begin{array}[]{ll}
E_{\ell (10)} \=  \frac{1}{2} \left(E_{\ell|1} - i E_{\ell|2}\right), &
E_{\ell(01)} \= \overline{E_{\ell(10)}},\\[6pt]
e_{-\ell(10)} \=  \frac{1}{2} \left(e_{-\ell|1} - i e_{-\ell|2}\right)
 , &e_{-\ell(01)} \= \overline{e_{-\ell(10)}},
 \end{array} \quad \ell = 0,1.\eeq
This is  the collection of  complex vector fields  that corresponds  to   the elements of the standard CR basis $\cB^{CR}$ by means of 
\eqref{absoluteparall}. From now,  we will   use the  notation $e_A$ and  $E_B$ to indicate just these  
vector fields. \par
\smallskip
The   vector fields $e_A$,  $E_B$ uniquely determine their dual  
 (real and complex) $1$-forms $\vartheta^A, \o^B$, defined by 
\beq \label{theta-omega} 
\vartheta^A(e_C) = \d^A_C,\quad \vartheta^A(E_D) = 0,\quad \o^B(e_A) = 0,\quad \o^B(E_D) = \d^B_D.
\eeq
Note that the $\gg$-valued $1$-form  $\varpi$   can be written in terms of  such $1$-forms as  
\beq\label{varpi} \varpi = \sum_A e^o_A \otimes \vartheta^A + \sum_B E^o_B \otimes \o^B.\eeq
\par
The vector fields $(e_A, E_B)$ and the dual $1$-forms $(\q^A, \o^B)$  have several geometric features, 
which derive from the special step-by-step construction of the Cartan connection $\varpi$ given in \cite{MS}. 
Let us  briefly recall them. \par
 First of all, we remind that  the  girdled CR manifold $(M, \cD, J)$ is  naturally equipped with a $J$-invariant, 2-dimensional, involutive  subdistribution $\cE$ of the distribution $\cD$, 
  defined at each point $x \in M$  by (see e.g. \cite{MS}, \S 2.1): 
\begin{multline} \cE_x\= \big\{v \in \cD_x:\text{there is vector field}\ X\ \text{in}\  \cD\ \\
\text{such that}\ X_x = v\ \text{and}\ \ [X, Y]_x \in \cD_x\ \text{for all vector fields $Y$ in}\ \cD\big\}.
\end{multline}
In other words, $\cE$ is the $J$-invariant distribution of  vector spaces, generated by the real 
 vectors  that are  in the kernels of the Levi forms of $(\cD, J)$.  \par
 In \cite{MS},  the bundle  $\pi: Q\to M$  is  obtained as the last step  of a tower of  three principal bundles,  one defined over the other,  as in the diagram
 $$Q = P^2 \overset{\pi^2} \longrightarrow P^1 \overset{\pi^1} \longrightarrow P^0\overset{\pi^0} \longrightarrow M,\qquad\text{with}\  \pi \= \pi^0 \circ \pi^1 \circ \pi^2.$$
  In turn,  each  bundle $P^i$ is  defined as a quotient $P^i = P^i_\sharp/N^i_\sharp $ by the action of a special group of matrices $N^i_\sharp$, of an appropriate principal 
 bundle $P^i_\sharp $ of linear frames of  the lower order bundle  \par
\hskip - 1cm{ \begin{picture}(300,110)(0,0)
 \setlength{\unitlength}{1pt}
 \linethickness{0.5mm}
 \put(282, 28){$M$}
\put(237, 68){$P^0_\sharp$}
\put(220, 94){$\cF r(M)$}
\put(237, 82){$\cup$}
 \linethickness{0.25mm}
\put(250, 30){\vector(1,0){30}}
 \linethickness{0.5mm}
 \put(250, 68){\vector(1,-1){30}}
\put(238, 19){\tiny $||$}
\put(237, 28){$P^0$}
\put(230, 4){$P^0_\sharp/N_\sharp$}
\put(185, 68){$P_\sharp^1$}
 \put(197, 68){\vector(1,-1){30}}
\put(168, 94){\small $\cF r(P^0)$}
\put(185, 82){$\cup$}
\put(185, 19){\tiny $||$}
\put(185, 28){$P^1$}
\put(175, 4){$P^1_\sharp/N^1_\sharp$}
 \linethickness{0.25mm}
\put(199, 30){\vector(1,0){30}}
\linethickness{0.5mm}
\put(130, 68){$P_\sharp^2$}
\put(112, 94){\small $\cF r(P^1)$}
\put(130, 82){$\cup$}
\put(132, 19){\tiny $||$}
\put(130, 28){$P^2$}
\put(120, 4){$P^2_\sharp/N^2_\sharp$}
 \put(147, 68){\vector(1,-1){30}}
  \linethickness{0.25mm}
 \put(145, 30){\vector(1,0){30}}
 \linethickness{0.5mm}
 \end{picture}}
 \par
\noindent The absolute parallelism $(e_A, E_B)$ on $Q = P^2$  is defined as  the unique frame field that takes values in a very 
special trivial subbundle $P^3_\sharp$ of the linear frame bundle $\cF r(P^2)$ of $P^2$.   \par
\smallskip
Each bundle of linear frames $P^i_\sharp \subset \cF r(P^{i-1})$,  $0 \leq i \leq 3$ (here, we set $M = P^{-1}$),  is determined by all linear frames of $P^{i-1}$
that are adapted to the natural distributions of $P^{i-1}$ (for instance, when $P^{i-1} = M$, the frames are adapted to the $J$-invariant   distributions  $\cE$ and $\cD$) and 
satisfy three  sets of conditions: 
\begin{itemize}
\item[a)] if the base point of the frame is  $u = [(f_i)] \in P^{i-1} = P^{i-1}_\sharp/N^{i-1}$, the first vectors of a frame in $P^i_\sharp|_u$ 
 are constrained  to  project onto one of the  linear frames $(f_i)$  in  $P^{i-1}_\sharp$, which  belong to the equivalence class $u = [(f_i)]$; 
\item[b)] the other vectors of a frame in $P^i_\sharp|_u$ must be  vertical with respect to the projection $\pi^i: P^{i-1} \to P^{i-2}$; 
\item[c)] any  linear frame in  $P^i_\sharp|_u$ is constrained by an appropriate set of normalising conditions; such conditions depend on which bundle $P^i_\sharp$ of linear frames we are considering -- we refer to \cite{MS} for the explicit formulation of such normalising conditions.  
 \end{itemize}
 Due  to (a), the frames in $P^i_\sharp $ not only satisfy the normalising conditions quoted in (c), but also  all conditions that 
 are residuals of the three types of  conditions  for  the frames in $P^{i-1}_\sharp$, in $P^{i -2}_\sharp$,  etc. In particular, the frame field   in  $P^3_\sharp$ that 
 gives the absolute parallelism $(e_A, E_B)$  on $Q = P^2$ satisfies a   set of conditions that inherits  from  the three types of constraints
 on the linear frames of the previous  steps.
Amongst such conditions one has that 
\begin{itemize}
\item[1)]  (the real and imaginary parts of) the vector fields $E_A$   are the infinitesimal transformations associated with (the real and imaginary parts of) the elements $E^o_A$,  determined by the right action 
of $H$ on $Q$; in particular, 
they are 
 generators of the  {\it vertical  distribution} $\cV \subset  TQ$, i.e.   the distribution of  the tangent spaces of the fibres of $\pi: Q \to M$; 
\item[2)] the distribution $\cH \subset TQ$, generated by  (the real and imaginary parts of) the vector fields  $e_A$, is such that for any $u \in Q$
the projection  $\pi_*|_u: T_u Q \to T_x M$, $x = \pi(u)$,  gives  a linear isomorphism 
 $\pi_*: \cH_u \to T_x M$ between $\cH_u$ and $T_x M$; 
\item[3)] for any $u \in Q$, 
 the   complex subspaces of $\cH_u^{\bC}$ 
 \begin{eqnarray} \label{d10}  \cD^{10 (\cH)}_u  \= \langle e_{-1(10)}|_u, e_{0(10)}|_u\rangle,& \cE^{10 (\cH)}_u  \= \langle e_{0(10)}|_u\rangle,\\
  \label{e10}  \cD^{01 (\cH)}_u \= \langle e_{-1(01)}|_u, e_{0(01)}|_u\rangle, & \cE^{01 (\cH)}_u \= \langle  e_{0(01)}|_u\rangle,
  \end{eqnarray}
project   isomorphically onto the  holomorphic spaces  $\cD^{10}_x\subset \cD^\bC_x$, $\cE^{10}_x\subset \cE^\bC_x$, $x = \pi(u)$, 
and  the  antiholomorphic spaces $\cD^{01}_x = \overline{\cD^{10}_x}$, $\cE^{01}_x = \overline{\cE^{10}_x}$, respectively.
  \end{itemize}
The other conditions  correspond to constraints on the curvature $2$-form $\k$ of the Cartan connection and will be  discussed in the next section.\par
\medskip
\subsection{The curvature constraints on the  Cartan connection} \label{constr}
  Consider now  the {\it curvature $2$-form} $\k$ of the Cartan connection $(Q, \varpi)$, that is  the $\so_{3,2}$-valued 2-form on $Q$,  defined by 
$$\k\= d \varpi + \textstyle\frac{1}{2} [\varpi, \varpi].$$
From basic properties of  Cartan connections and the fact that the vector fields $E_A$ are infinitesimal transformations on $Q$, corresponding  to the elements $E^o_A \in \gh = Lie(H)$,  
the expansion of $\k$ in terms of the  pointwise linearly independent  2-forms $(\q^A \wedge \q^B, \q^A \wedge \o^C$, $\o^C \wedge \o^D)$,  determined by  the dual coframe  \eqref{theta-omega},  
has necessarily the form
\beq \label{kappa} \k = d \varpi +\textstyle \frac{1}{2}  [\varpi, \varpi] =  \sum T^A_{BC} e^o_A \otimes \q^B \wedge \q^C + \sum R^{D}_{BC} E^o_{D} \otimes \q^B \wedge \q^C,\eeq
for appropriate (real and complex) functions $T^A_{BC}$ and $R^{D}_{BC}$. \par
\smallskip
The curvature components $T^A_{BC}$ and $R^{D}_{BC}$  are determined by the   Lie brackets of pairs of  vector fields $e_A$ as follows.  Let us denote by 
 $(\cdot)^{e^o_A}: \so_{3,2} \to \langle e^o_A\rangle$ and $(\cdot)^{E^o_B}: \so_{3,2} \longrightarrow \langle e^o_B\rangle$
 the standard projections of $\so_{3,2}$ along the vectors of the basis $\cB^{CR}$ and by $ \gc^A_{BC} {\=} \ ([e^o_B, e^o_C])^{e^o_A}$ and $\gd^D_{BC} {\=} \ ([e^o_B, e^o_C])^{E^o_D}$ the structure constants of $\so_{3,2}$ in the basis  $\cB^{CR}$. Then,   by Koszul formula for exterior derivatives and the definition of $\k$,  we   have  
\beq \label{T}
 T^A_{BC} = d \q^A(e_B, e_C) +   ([e^o_B, e^o_C])^{e^o_A} = - \q^A([e_B, e_C]) + \gc^A_{BC},\eeq
\beq  \label{R} R^{D}_{BC} = d \o^D(e_B, e_C) +   ([e^o_B, e^o_C])^{E^o_D} =  - \o^{D}([e_B, e_C]) + \gd^{D}_{BC}.\eeq
By comparison with \eqref{alternative}, we immediately see that,  modulo the structures constants of $\so_{3,2}$,  {\it the curvature components  $T^A_{BC}$ and $R^{D}_{BC} $ of the curvature  $2$-form $\k$ are nothing but the  
structure  functions   of the absolute parallelism $(e_A, E_B)$}.  In fact, this   is  a well known general fact  on Cartan connections.\par  
\smallskip
As mentioned above, besides the conditions (1)  -- (3) of \S \ref{sect3.3}, the absolute parallelism $(e_A, E_B)$ is constrained by other  normalising conditions. 
They are conditions on the Lie brackets between the vectors $e_A$ and, through \eqref{T} and \eqref{R}, they can all be expressed in terms of the 
curvature components $T^A_{BC}$ and $R^D_{BC}$.  For readers convenience, we give  
here  the complete  list of such   constraints  and we refer to   \cite{MS} for further details. 
  \par
 \medskip
 \noindent{\it Integrability of the complex structure and involutivity of $\cE^\bC$}.\par
 \noindent  
From  \eqref{T} and the fact that  the  complex distributions   $\cE^{10(\cH)} + \cE^{01(\cH)}$ and $\cD^{10(\cH)}$,   defined in \eqref{d10} and \eqref{e10},   project  onto the involutive distributions  $\cE^\bC$ and $\cD^{10}$ of $M$, 
one has 
\beq \label{primoconstr} 
\begin{array}{l}T^{A}_{0(10)\, 0(01)}  = 0 \  \text{for}\ A \in \{-2, -1(10), -1(01)\},\\
T^{A'}_{i(10) j (10)}  = \overline{T^{A'}_{i(10) j (10)}} = 0\  \text{for}\  i, j \in \{-1, 0\}, \ A' \in \{-2, -1(01), 0(01)\} .
\end{array}
\eeq
\par 
\smallskip
 \noindent{\it The  distribution $\cE^{10}$  is in the kernel of  Levi forms}.\par
 \noindent  
From \eqref{T}  and the fact that the spaces   $\cE^{10(\cH)}_u$, $\cE^{01(\cH)}_u$  project into the  kernel 
of the Levi form,   one has 
\beq
 T^{-2}_{-1(01)\,0(10)} = T^{-2}_{-1(10)\,0(01)} = 0 .
\eeq
\par
\smallskip
 \noindent{\it Normalising conditions on the  frames in  $P^0_\sharp$}.\par
 \noindent 
From \eqref{T} and condition (5.2) in \cite{MS}, one has  
 \beq T^{-2}_{-1(10)\,-1(01)} = T^{-1(10)}_{-1(01)\,0(10)} = T^{-1(01)}_{-1(10)\,0(01)}  = 0.\eeq
 \par
\smallskip
 \noindent{\it Normalising conditions on the  frames of $P^1_\sharp$}.\par
 \nopagebreak
 \noindent  
 From \eqref{T} and the normalising conditions in \cite{MS}, given by formula (6.7) and  condition $\b_K = 0$ after Lemma 6.5 of that paper, one has 
\beq
 \begin{split}
 & T^{-1(10)}_{-1(10)\, 0(10)} = T^{-1(10)}_{-1(10)\, 0(01)} = T^{-1(01)}_{-1(01)\, 0(01)} = T^{-1(01)}_{-1(01)\, 0(10)} = 0,\\[2pt]
&T^{0(10)}_{-1(01)\, 0(01)} = T^{0(01)}_{-1(10)\, 0(10)} = 0.
\end{split}
\eeq
  \par
\smallskip
  \noindent{\it Property of the strongly adapted frames in $P^1_\sharp$}.\par
 \noindent  
From \eqref{T} and Lemma 6.6 (ii) in \cite{MS} one has  
 \beq T^{-2}_{-2\, 0(10)} = T^{-2}_{-2\, 0(01)} = 0.\eeq
 \par
\smallskip
  \noindent{\it Normalising conditions on the strongly adapted frames in $P^2_\sharp$}.\par
 \noindent  
 From \eqref{T} and  the normalising conditions (7.2)  on $\g_K$ in \cite{MS}, one has
 \beq T^{0(10)}_{-1(10)\, 0(10)} = T^{0(01)}_{-1(01)\, 0(01)} = T^{0(10)}_{-1(10)\, 0(01)} = T^{0(01)}_{-1(01)\, 0(10)} = 0.\eeq
  \par
\smallskip
  \noindent{\it Normalising conditions on the strongly adapted frames in $P^3_\sharp$}.\par
 \noindent 
 From \eqref{R} and the normalising conditions (8.2)  on $\varepsilon_K$ in \cite{MS}, one has 
 \beq  \label{ultimoconstr} R^{0(10)}_{-2\, 0(10)} =  R^{0(01)}_{-2\, 0(01)} = R^{0(01)}_{-2\, 0(10)} =  R^{0(10)}_{-2\, 0(01)} = 0.\eeq
 \par
 \smallskip
 Besides  \eqref{primoconstr} -- \eqref{ultimoconstr},  the absolute parallelism is subjected to    three further conditions of cohomological nature. They are
 \begin{itemize}
 \item[a)] the condition given  in (6.21) of \cite{MS}, which  is equivalent  to a  system of linear equations on 
 $T^{-2}_{-2\, -1(10)}$, $ T^{-2}_{-2\, -1(01)}$ and $T^{-1(10)}_{-1(10)\,\ -1(01)}$;
 \item[b)] the condition given in (7.4) in \cite{MS}, which   is equivalent  to a  system of linear equations on 
 $T^{-1(10)}_{-2\, -1(10)}$, $T^{-1(10)}_{-2\, -1(01)}$, $T^{0(10)}_{-1(10)\,-1(01)}$,  $R^{0(10)}_{-1(10)\, -1(01)}$
 and their complex conjugates; 
\item[c)] the condition  given in (8.4) in \cite{MS}, which  is equivalent  to a  system of linear equations on 
$ T^{0(10)}_{-2\, -1(10)}$, $T^{0(10)}_{-2\, -1(01)}$, 
 $R^{0(10)}_{-2(10)\, -1(01)}$, $R^{0(10)}_{-2\, -1(01)}$,  $R^{1(10)}_{-1(10)\,-1(01)}$
 and their complex conjugates.
 \end{itemize}
 The  explicit expressions for the linear systems corresponding to the constraints  (a), (b), (c) can be determined with  straightforward computations. 
 An exposition of  such computations, which uses only elementary tools, can be found in Appendix \ref{appendixB}. 
 The result  is that the constraints (a), (b) and (c) are equivalent to  the linear equations 
\begin{align} & \text{a)}\qquad T^{-2}_{-2\, -1(10)} =  T^{-2}_{-2\, -1(01)} =  T^{-1(10)}_{-1(10)\,-1(01)} = 0,\hskip56pt\ \\
 \nonumber &\text{b)}\qquad  T^{-1(10)}_{-2\, -1(10)} =  T^{-1(10)}_{-2\, -1(01)} =T^{0(10)}_{-1(10)\, -1(01)} =T^{-1(01)}_{-2\, -1(01)} = T^{-1(01)}_{-2\, -1(10)} = \\
& \hskip 2cm  
 =T^{0(01)}_{-1(10)\, -1(01)}  = R^{0(10)}_{-1(10)\,-1(01)} =  R^{0(01)}_{-1(10)\,-1(01)} = 0,\\
\nonumber & \text{c)}\qquad   \overline{R^{0(10)}_{-2\,-1(10)} }   =   - \frac{1}{2} T^{0(10)}_{-2\,-1(10)}   - \frac{1}{2}   R^{0(10)}_{-2\,-1(01)}, \\
& \label{c} \hskip 1.2cm R^{1(10)}_{-1(10) \,-1(01)} =  \frac{i}{2} T^{0(10)}_{-2\,-1(10)}  - \frac{i}{2} R^{0(10)}_{-2 \,-1(01)}
\end{align}
and to the equations that  follows from \eqref{c}  by complex conjugation.
\par
\subsection{The structure equations of a girdled CR manifold}
The projections of the values of curvature $\k$ along each element  of the basis $\cB^{CR}$ 
give explicit expressions 
for   the exterior differentials $d\q^A$  and $d\o^B$ in terms  of the  pointwise linearly independent  2-forms $(\q^A \wedge \q^B, \q^A \wedge \o^C$, $\o^C \wedge \o^B)$, 
i.e. the {\it structure equations of the absolute parallelism $(e_A, E_B)$} (see \S \ref{sect2.2}). Here is the complete list of 
these structure equations, where we set equal to $0$ all  terms  $T^A_{BC}$ that  are bound  to  vanish by  the curvature  constraints  in \S \ref{constr}.\par
\beq  \label{str1} d \vartheta^{-2} + \frac{i}{2} \vartheta^{-1(10)} \wedge \vartheta^{-1(01)} - \left(\o^{0(10)} + \o^{0(01)}\right)\wedge \vartheta^{-2}= 0,
\eeq
 \begin{multline}  \label{str2}  d \vartheta^{-1(10)} - \vartheta^{0(10)}\wedge \vartheta^{-1(01)} - \o^{0(10)}  \wedge \vartheta^{-1(10)} +  i \o^{1(10)}\wedge \vartheta^{-2} = 
 \\[1pt]
 = T^{-1(10)}_{-2 \,0(10)} \vartheta^{-2} \wedge \vartheta^{0(10)} +
T^{-1(10)}_{-2 \,0(01)} \vartheta^{-2} \wedge \vartheta^{0(01)}, \end{multline}
 \begin{multline}  d \vartheta^{0(10)} - \left(\o^{0(10)} -\o^{0(01)}\right)\wedge \vartheta^{0(10)} + \frac{1}{2} \o^{1(10)} \wedge \vartheta^{-1(10)} = 
 \\[1pt]
 = T^{0(10)}_{-2 \,-1(10)} \vartheta^{-2} \wedge \vartheta^{-1(10)} +  T^{0(10)}_{-2 \,-1(01)} \vartheta^{-2} \wedge \vartheta^{-1(01)} +
\\[1pt]  \label{str3}
+  T^{0(10)}_{-2 \,0(10)} \vartheta^{-2} \wedge \vartheta^{0(10)} + T^{0(10)}_{-2 \,0(01)} \vartheta^{-2} \wedge \vartheta^{0(01)} +T^{0(10)}_{-1(01) \,0(10)} \vartheta^{-1(01)} \wedge \vartheta^{0(10)} + \\[1pt]
  +T^{0(10)}_{0(10) \,0(01)} \vartheta^{0(10)} \wedge \vartheta^{0(01)},
 \end{multline}
 \begin{multline}  d \o^{0(10)} - \vartheta^{0(10)} \wedge \vartheta^{0(01)}+ \frac{1}{2} \o^{1(01)} \wedge \vartheta^{-1(10)} + 
\o^2 \wedge \vartheta^{-2} = 
 \\[1pt]
 =\underset{\text{\bf constrained by \eqref{c}}}{R^{0(10)}_{-2 \,-1(10)} }\vartheta^{-2} \wedge \vartheta^{-1(10)} +  R^{0(10)}_{-2 \,-1(01)} \vartheta^{-2} \wedge \vartheta^{-1(01)}+\\[1pt]
 \label{str4}
  + R^{0(10)}_{-1(10) \,0(10)} \vartheta^{-1(10)} \wedge \vartheta^{0(10)} +
   R^{0(10)}_{-1(10) \,0(01)} \vartheta^{-1(10)} \wedge \vartheta^{0(01)}+\\[1pt]
     +R^{0(10)}_{-1(01) \,0(10)} \vartheta^{-1(01)} \wedge \vartheta^{0(10)} + R^{0(10)}_{-1(01) \,0(01)} \vartheta^{-1(01)} \wedge \vartheta^{0(01)}+ \\[1pt]
  +R^{0(10)}_{0(10) \,0(01)} \vartheta^{0(10)} \wedge \vartheta^{0(01)},
 \end{multline}
 \begin{multline}  \nonumber d \o^{1(10)} - \o^{1(01)}\wedge \vartheta^{0(10)} - \o^{1(10)}\wedge\o^{0(01)}  + i \o^2 \wedge \vartheta^{-1(10)} =  \\[1pt]
 =R^{1(10)}_{-2 \,-1(10)} \vartheta^{-2} \wedge \vartheta^{-1(10)} +  R^{1(10)}_{-2 \,-1(01)} \vartheta^{-2} \wedge \vartheta^{-1(01)} +\\[1pt]
 \nonumber   +  R^{1(10)}_{-2 \,0(10)} \vartheta^{-2} \wedge \vartheta^{0(10)} + R^{1(10)}_{-2 \,0(01)} \vartheta^{-2} \wedge \vartheta^{0(01)} +\\[1pt]
  +\underset{\text{\bf constrained by \eqref{c}}}{ R^{1(10)}_{-1(10) \,-1(01)}} \vartheta^{-1(10)} \wedge \vartheta^{-1(01)} + \\[1pt]
  + R^{1(10)}_{-1(10) \,0(10)} \vartheta^{-1(10)} \wedge \vartheta^{0(10)} +
   R^{1(10)}_{-1(10) \,0(01)} \vartheta^{-1(10)} \wedge \vartheta^{0(01)}+
   \end{multline}
 \begin{multline}\label{str5} 
    +R^{1(10)}_{-1(01) \,0(10)} \vartheta^{-1(01)} \wedge \vartheta^{0(10)} + R^{1(10)}_{-1(01) \,0(01)} \vartheta^{-1(01)} \wedge \vartheta^{0(01)}+ \\[1pt]
  +R^{1(10)}_{0(10) \,0(01)} \vartheta^{0(10)} \wedge \vartheta^{0(01)},
 \end{multline}
 \begin{multline} \nonumber d \o^2 -  \frac{i}{2} \o^{1(10)} \wedge \o^{1(01)} + \left(\o^{0(10)} +  \o^{0(01)}\right) \wedge  \o^2 =  \\[1pt]
 =R^{2}_{-2 \,-1(10)} \vartheta^{-2} \wedge \vartheta^{-1(10)} +  R^{2}_{-2 \,-1(01)} \vartheta^{-2} \wedge \vartheta^{-1(01)} +
\\[1pt]
+  R^{2}_{-2 \,0(10)} \vartheta^{-2} \wedge \vartheta^{0(10)} + R^{2}_{-2 \,0(01)} \vartheta^{-2} \wedge \vartheta^{0(01)} +\\[1pt]
  +R^{2}_{-1(10) \,-1(01)} \vartheta^{-1(10)} \wedge \vartheta^{-1(01)} + 
\end{multline}
\begin{multline}
 \label{str6}
  + R^{2}_{-1(10) \,0(10)} \vartheta^{-1(10)} \wedge \vartheta^{0(10)} +
   R^{2}_{-1(10) \,0(01)} \vartheta^{-1(10)} \wedge \vartheta^{0(01)}+\\[1pt]
    +R^{2}_{-1(01) \,0(10)} \vartheta^{-1(01)} \wedge \vartheta^{0(10)} + R^{2}_{-1(01) \,0(01)} \vartheta^{-1(01)} \wedge \vartheta^{0(01)}+ \\[1pt]
  +R^{2}_{0(10) \,0(01)} \vartheta^{0(10)} \wedge \vartheta^{0(01)}.
 \end{multline}
 \par
\medskip
\subsection{Comparison with  other absolute  parallelisms}
As we already mentioned, other canonical absolute parallelisms 
for girdled CR manifolds, not associated with  Cartan connections, have been recently  given in \cite{IZ, Po}.  Note  also that the absolute  parallelism   in \cite{Po} is 
defined only for the  girdled CR manifolds admitting   no local equivalence  with 
the homogeneous girdled CR manifold $M_o$. \par
Let us now focus on the  canonical absolute parallelism $(P, (X_i), \wt{(\cdot)})$,  defined  in  \cite{IZ} for an arbitrary girdled CR manifold $(M, \cD, J)$.
There, the  bundle $\pi: P \to M$ has
 $5$-dimensional fibers, but it has no natural structure of principal bundle over $M$.  
 The  absolute parallelism $(X_i)_{i = 1}^{10}$ on $P$ is associated with a dual  coframes field,   given by the real and imaginary parts of 
ten   $\bC$-valued $1$-forms,  denoted by 
$(\o, \o^1, \overline{\o^1}, \varphi^2, \overline{\f^2},  \theta^2, \overline{\theta^2},  \varphi^1, \overline{\varphi^1}, \psi)$ and with 
$\o$ and $\psi$  taking only imaginary values. \par 
Since the bundle $P$  and  the principal bundle $Q$ of our Cartan connection  have  the same dimension,   if we consider them  as mere 
 bundles  with no further structures, we may locally identify them.  From the construction of $P$, we may also assume that, 
 under this identification,  the $1$-form $\o$ is equal to    $\o  = -2  i \q^{-2}$, where $\q^{-2}$ is the  $1$-form of our parallelism, 
 defined in \S \ref{sect3.3}.  \par
We now recall that the $1$-forms of the absolute parallelism in \cite{IZ} are characterised by  the fact that they satisfy an appropriate set of structure equations. 
The first two of this set are 
\begin{align} \label{11} d \o &= - \o^1 \wedge \o^{\bar 1} - \o \wedge (\varphi^2 + \varphi^{\bar 2}),\\
 \label{12}  d \o^1 &= \theta^2 \wedge \o^{\bar 1} - \o^1 \wedge \f^2 - \o \wedge \f^1.
 \end{align}
Comparing them with  our structure equations  \eqref{str1} and \eqref{str2}  and through a tedious but straightforward computation, one can  check  that   
 the equations \eqref{11} and \eqref{12} are satisfied by the  $1$-forms on $P \simeq Q$, defined by 
\beq 
\label{3.37}
 \begin{split}
\o &  \= - 2 i  \q^{-2},\\
\o^1 &  \=  \q^{-1(10)}  - \overline{T^{-1(10)}_{-2\,0(10)}}\q^{-2},\\
     \f^2 &\= \o^{0(10)} + \frac{i}{2}\overline{T^{-1(10)}_{-2\, 0(10)}}  \q^{-1(01)} ,\\
\theta^2 & \= \q^{0(10)} + i \overline{T^{-1(10)}_{-2\,0(10)}}\left( \q^{-1(10)}  - \overline{T^{-1(10)}_{-2\,0(10)}}\q^{-2}\right),
\\
\f^1 & \=  \frac{1}{2} \o^{1(10)}Ê - \frac{i}{2}
T^{-1(10)}_{-2 \,0(01)}  \vartheta^{0(01)}-  \frac{1}{2}\overline{T^{-1(10)}_{-2\,0(10)}}  T^{-1(10)}_{-2\,0(10)}  \q^{-1(10)}   + \\
&  +  \frac{1}{4}\overline{T^{-1(10)}_{-2\,0(10)}} \overline{T^{-1(10)}_{-2\,0(10)}} \q^{-1(01)}  - \frac{i}{2}  \overline{T^{-1(10)}_{-2\,0(10)}}\o^{0(01)}  -  \frac{i}{2} d\overline{T^{-1(10)}_{-2\,0(10)}}.
\end{split}
\eeq
Now,  we expect that if the $1$-forms \eqref{3.37} are appropriately modified with additional terms involving $\o^{1(10)}, \o^{1(01)}$ and $\o^2$, they will satisfy not only the first two structure equations of the absolute parallelism in \cite{IZ}, but  also  all other  structure equations
of that parallelism. \par
On the base of this expectation, the construction in \cite{IZ}  seems to start diverging   from ours 
precisely when  the absolute parallelism is required to satisfy  \eqref{12}.  In fact, this is a constraint that 
amounts to imposeÊ that the curvature components  $T^{-1(10)}_{-2\,0(10)}$ and $T^{-1(10)}_{-2\,0(01)}$ are 
absorbed into the definition of the vector fields of the absolute parallelism. Since these curvature  components are not 
invariant under the right action of the structure group  of  $\pi: Q \to M$,   the  Êconstraint given by \eqref{12}  is  plausibly one of the  main 
reasons for the fact that the constructive process  in \cite{IZ} does not produce    a Cartan connection. \par
An analogous comparison  between  our canonical Cartan connection and the parallelism in \cite{Po}  might be done following the same   line of arguments. 
We leave this  task to the interested reader.\par
\appendix
\section{} \label{appendixB}
\setcounter{equation}{0}
\subsection{The Cartan-Killing form of $\so_{3,2}$}
For the  following computations, it turns out that the standard basis $\cB^o$ of $\so_{3,2}$, defined in \eqref{basisso32}, is not very convenient. 
In place of that basis, it is by far more useful to consider a new basis $\cB = (f_\a)_{\a = 1, \ldots, 10}$, with elements   
\beq
{\scalefont{0.88}
\begin{aligned}\label{basisR}
f_1 & \= \frac{1}{\sqrt{6}} e_{-2}, \quad f_2  \= \frac{1}{\sqrt{6}}(e_{-1(10)}  + e_{-1(01)}),  \quad f_3 \= \frac{i}{\sqrt{6}} (e_{-1(10)} - e_{-1(01)}), \\
  f_4 & \= \frac{1}{\sqrt{12}}(e_{0(10)} + e_{0(01)}), \quad  f_5 \=  \frac{i}{\sqrt{12}} (e_{0(10)} -e_{0(01)}),\\
f_6 &\=  \frac{1}{\sqrt{12}} (E_{0(10)} + E_{0(01)}),  \quad   f_7 \= \frac{i}{\sqrt{12}} ( E_{0(10)} -  E_{0(01)}), \\ 
 f_8 & \= \frac{1}{\sqrt{6}} (E_{1(10)} + E_{1(01)}), \quad  f_9 \=  \frac{i}{\sqrt{6}} ( E_{1(10)} - E_{1(01)}), \quad
 f_{10}  \= \frac{1}{\sqrt{6}} E^2 .
\end{aligned}
}
\eeq
The main motivation for considering such new basis comes from the fact  the entries of    the Cartan-Killing form $\langle \cdot , \cdot \rangle$ of $\so_{3,2}$ in this basis  are equal   to $\pm 1$ or $0$. More precisely, using Table 1 in \cite{MS}, one can  check that  
the components of  $\langle \cdot , \cdot \rangle$  in the basis  $\cB$ are 
\beq  \label{CartanMatrix} {\scalefont{0.75}\arraycolsep=3pt
\left(\begin{array}{cccccccccc}
0 &  0 &  0 &  0 &  0 &  0 &  0 &  0 &  0 &  1\\
0 &  0 &  0 &  0 &  0 &  0 &  0 &  1 &  0 &  0\\
0 &  0 &  0 &  0 &  0 &  0 &  0 &  0 &  1 &  0\\
0 &  0 &  0 &  1 &  0 &  0 &  0 &  0 &  0 &  0\\
0 &  0 &  0 &  0 &  1 &  0 &  0 &  0 &  0 &  0\\
0 &  0 &  0 &  0 &  0 &  1 &  0 &  0 &  0 &  0\\
0 &  0 &  0 &  0 &  0 &  0 &  -1 &  0 &  0 &  0\\
0 &  1 &  0 &  0 &  0 &  0 &  0 &  0 &  0 &  0\\
0 &  0 &  1 &  0 &  0 &  0 &  0 &  0 &  0 &  0\\
1 &  0 &  0 &  0 &  0 &  0 &  0 &  0 &  0 &  0
\end{array}\right).}\eeq
\subsection{The space $\ker \p^*|_{C^2_1(\gm_, \gg)}$}
\label{sectA2}
We now want  to  show that  the action of the codifferential $\p^*$ of $\so_{3,2}$ on the bilinear maps  of shifting degree $+1$  in 
$\Hom(\L^2 \gm_{-} , \gg)$, $\gm_- \= \gg_{-2} + \gg_{-1}$,  has trivial kernel.  As above,
 we denote by $\langle \cdot , \cdot \rangle$ the Cartan-Killing form of $\so_{3,2}$, by $\cB = (f_\a)$ the basis defined in \eqref{basisR} and by $\cB^*  = (f^\a)$ its dual basis. 
Finally, for each element $f_\a \in \cBì$, we denote by  $\wh f_\a$ the unique element  in $\cB$  such that $f^{\a} = \pm \langle \wh f_\a, \cdot\rangle$ and by $\wh f^\a$ the corresponding dual element in $\cB^*$.  The rest of the notation is taken from \cite{MS}.  \par
Consider a  bilinear map $\t \in \Hom(\L^2 \gm_-, \gg)$  of shifting degree $+1$:
\beq \label{4.1}\t = \t^1_{12} f_1 \otimes (f^1 \wedge f^2) + \t^1_{13} f_1 \otimes (f^1 \wedge f^3) + \t^2_{23} f_2 \otimes (f^2 \wedge f^3) +  \t^3_{23} f_3 \otimes (f^2 \wedge f^3) .\eeq
By definition of $\p^*$, this tensor is in $\ker \p^*$ if and only if 
\beq \label{4.2} \langle \p^* \t, A \rangle = -  \langle \t, \p A \rangle = 0\eeq
for any  $A = A^\a_\b f_\a \otimes f^\b \in \Hom(\gh, \gg)$. From \eqref{4.1},  equation  \eqref{4.2} is equivalent to a linear equation on the    $\t^i_{jk}$ whose non trivial  coefficients are
$$\wh f^1(\p A(\wh f_1, \wh f_2)) ,\quad \wh f^1( \p A(\wh f_1, \wh f_3)),\quad  \wh f^2(\p A(\wh f_2, \wh f_3)),\quad  \wh f^3(\p A(\wh f_2, \wh f_3)).$$
The computation of $\wh f^1(\p A(\wh f_1, \wh f_2))$ is straightforward and gives 
\begin{align*}
\wh f^1(\p A(\wh f_1, \wh f_2)) &= f^{10}(\p A(f_{10}, f_8)) =\\
& = f^{10}\bigg( f_{10} \cdot A(f_8) - f_8 \cdot A(f_{10}) - A([f_{10}, f_8])\bigg) = \\
& =   A^\a_8 f^{10}(\ad_{f_{10}}(f_\a)) - A^\a_{10} f^{10}(\ad_{f_{8}}(f_\a)) - f^{10}(A([f_{10}, f_{8}])) = \\
& =  -  \frac{1}{\sqrt{3}} A^6_8  + \frac{1}{\sqrt{6}} A^9_{10}.
\end{align*}
Similar computations give all other coefficients of the equations and 
 \eqref{4.2} reduces to
\beq 
{\scalefont{0.90}
\begin{aligned}
& \t^1_{12}\big(-  \frac{1}{\sqrt{3}} A^6_8  + \frac{1}{\sqrt{6}} A^9_{10}\big) +  \t^1_{13}\big( -  \frac{1}{\sqrt{3}} A^6_9  - \frac{1}{\sqrt{6}} A^8_{10}\big) +\\
& +  \t^2_{23} \big(-  \frac{1}{2\sqrt{3}} A^6_9  -  \frac{1}{2\sqrt{3}} A^4_9 +  \frac{1}{2\sqrt{3}} A^7_8  +  \frac{1}{2\sqrt{3}} A^5_8  +\frac{1}{\sqrt{6}}A^8_{10}  \big) + \\
& + \t^3_{23}\big(  \frac{1}{2\sqrt{3}} A^7_9  -  \frac{1}{2\sqrt{3}} A^5_9 +  \frac{1}{2\sqrt{3}} A^6_8  -  \frac{1}{2\sqrt{3}} A^4_8
+\frac{1}{\sqrt{6}}A^9_{10}
\big) = 0.
\end{aligned}
}
\eeq
By arbitrariness of  $A$,   it follows that  $\t \in \ker \p^*$ if and only if 
$ \t^1_{12}   = \t^1_{13} = \t^2_{23} = \t^3_{23} = 0$,  meaning that  $\ker \p^*|_{C^2_1(\gm_, \gg)} = 0$, as claimed. 
\par
\medskip
\subsection{The space $(\p \gl^1)^\perp$} We recall that, according to Lemma 6.5 in \cite{MS} and the definition of the (abelian) group $L^1$, the abelian Lie algebra $\gl^1 = Lie(L^1)$ 
can be identified with the real vector space generated by the linear maps
\beq
{\scalefont{0.90}
\begin{aligned}
& B_1 \= e_{-1(10)} \otimes e^{-2} +  e_{-1(01)} \otimes e^{-2},\\
& B_2 \= i e_{-1(10)} \otimes e^{-2} - i  e_{-1(01)} \otimes e^{-2},\\
& B_3\=\left(e_{0(10)} - E_{0(01)}\right)\otimes e^{-1(10)} +  \left(e_{0(01)} - E_{0(10)}\right)\otimes e^{-1(01)},\\
& B_4 \= i \left(e_{0(10)} - E_{0(01)}\right)\otimes e^{-1(10)} - i \left(e_{0(01)} - E_{0(10)}\right)\otimes e^{-1(01)},\\
& B_5\= i \left( e_{0(10)} + E_{0(01)}\right)\otimes e^{-1(10)} - i \left( e_{0(01)} + E_{0(10)}\right)\otimes e^{-1(01)},\\
& B_6\=  \left( e_{0(10)} + E_{0(01)}\right)\otimes e^{-1(10)} - \left( e_{0(01)} + E_{0(10)}\right)\otimes e^{-1(01)},\\
& B_7\= \left(E_{0(10)} + E_{0(01)}\right)\otimes e^{-1(10)} + \left(E_{0(10)} + E_{0(01)}\right)\otimes e^{-1(01)},\\
& B_8\= i\left(E_{0(10)} + E_{0(01)}\right)\otimes e^{-1(10)} -i  \left(E_{0(10)} + E_{0(01)}\right)\otimes e^{-1(01)}.
\end{aligned}
}
\eeq
We also recall that the elements  $\t^1 \in Tor^1(\gm)$  have the form
form
$$\t^1 = \t^{-2}_{-2 -1(10)} e_{-2} \otimes e^{-2} \wedge e^{-1(10)} + \overline{\t^{-2}_{-2 -1(10)}} e_{-2} \otimes e^{-2} \wedge e^{-1(01)}+$$
$$ +  \t^{-1(10)}_{-1(10) -1(01)} e_{-1(10)} \otimes e^{-1(10)*} \wedge e^{-1(01)*} +$$
$$  + \overline{\t^{-1(01)}_{-1(10) -1(01)} }e_{-1(01)} \otimes e^{-1(10)*} \wedge e^{-1(01)*}.$$
Now, let us  choose as  $\ad_{E^0_2}$-invariant inner product on the space $Tor^1(\gm)$  the  sum of  an arbitrary inner product on $\gm_{-2}$ and the standard Hermitian inner product of $\bC \simeq \gm_{-1}$. 
This implies that, in order to determine the subspace  $(\p \gl^1)^\perp \subset  Tor^1(\gm)$, 
the only relevant components of the generators $\p B_i$ of $\p \gl^1$ are  the components 
$$(\partial B_i)^{-2}_{-2 -1(10)}  \=  e^{-2} \left(\p B_i(e_{-2}, e_{-1(10)})\right),$$
$$(\partial B_i)^{-1(10)}_{-1(10) -1(01)}   \=  e^{-1(10)} \left(\p B_i(e_{-1(10)}, e_{-1(01)})\right).$$
 We now observe that 
  \begin{align*}
  (\partial B_3)^{-2}_{-2 -1(10)} &= e^{-2} \left([e_{-2}, B_3(e_{ -1(10)})] - [e_{-1(10)}, B_3(e_{-2})]\right) = -1,\\
(\partial B_4)^{-2}_{-2 -1(10)} &= e^{-2} \left([e_{-2}, B_4(e_{ -1(10)})] - [e_{-1(10)}, B_4(e_{-2})]\right) = -i,
  \end{align*}
meaning  that $\p \gl^1$ contains the 2-dimensional real subspace generated  by 
$$ e_{-2} \otimes e^{-2} \wedge e^{-1(10)} + \overline{e_{-2} \otimes e^{-2} \wedge e^{-1(10)}},$$
 $$ i e_{-2} \otimes e^{-2} \wedge e^{-1(10)} -i \overline{e_{-2} \otimes e^{-2} \wedge e^{-1(10)}}.$$
 This yields that if $\t^1 \in (\p\gl^1)^\perp$, then $ \t^{-2}_{-2 -1(10)} = \overline{\t^{-2}_{-2 -1(10)}} = 0$. 
Similar computations show that 
$$(\partial B_1)^{-1(10)}_{-1(10) -1(01)} = - \frac{i}{2},\qquad 
(\partial B_2)^{-1(01)}_{-1(10) -1(01)} =  - \frac{1}{2},
$$
hence  that $\p \gl^1$ contains the 2-dimensional real subspace generated  by 
$$ e_{-1(10)} \otimes e^{-1(10)} \wedge e^{-1(01)} + \overline{e_{-1(10)} \otimes e^{-1(10)} \wedge e^{-1(01)}},$$
 $$ i e_{-1(10)} \otimes e^{-1(10)} \wedge e^{-1(01)} - i \overline{e_{-1(10)} \otimes e^{-1(10)} \wedge e^{-1(01)}},$$
and therefore that if $\t^1 \in (\p\gl^1)^\perp$, then $ \t^{-1(10)}_{-1(10) -1(01)} = \overline{ \t^{-1(10)}_{-1(10) -1(01)}} = 0$. 
We therefore conclude that 
$(\p\gl^1)^\perp  = 0$.\par Since $\ker \p^*|_{C^2_1(\gm_, \gg)}$ is  trivial as well (see  \S \ref{sectA2}), 
condition (6.21) in \cite{MS} is equivalent to requiring that  {\it the $c$-torsion $c^1_K$ is identically  equal to $0$}.
\par
\medskip
\subsection{The space $\ker \p^*|_{C^2_2(\gm_-, \gg)}$} Here, we want show that the space of the bilinear maps in $\Hom(\L^2 \gm_{-} , \gg)$ of shifting degree $+2$
that are  in $\ker \p^*|_{C^2_2(\gm_, \gg)}$, is trivial.   This amount to say  that condition (7.4) of \cite{MS} reduces to 
$c^2_K  = 0$.\par
\smallskip
A bilinear map  $\t \in \Hom(\L^2 \gm_-, \gg)$  of shifting degree $+2$ has  the form  
\beq  \label{6.1}
{\scalefont{0.88}
\begin{aligned}
\t &{=} \t^2_{12} f_2 \otimes (f^1 \wedge f^2) + \t^3_{12} f_3 \otimes (f^1 \wedge f^2) +\t^2_{13} f_2 \otimes (f^1 \wedge f^3) + \t^3_{13} f_3 \otimes (f^1 \wedge f^3) +\\
& +\t^4_{23} f_4 \otimes (f^2 \wedge f^3) +  \t^5_{23} f_5 \otimes (f^2 \wedge f^3) + \t^6_{23} f_6 \otimes (f^2 \wedge f^3) +  \t^7_{23} f_7 \otimes (f^2 \wedge f^3)
 .
 \end{aligned}
 }
 \eeq
As in \S \ref{sectA2}, this tensor is in $\ker \p^*$ if and only if 
$ \langle \p^* \t, A \rangle = -  \langle \t, \p A \rangle = 0$
for any element $A   = A^\a_\b f_\a \otimes f^\b\in \Hom(\gh, \gg)$. By \eqref{6.1},  this corresponds to a linear equation on  the components $\t^i_{jk}$ with  coefficients 
\beq \nonumber
{\scalefont{0.90}
\begin{aligned}
&\wh f^2(\p A(\wh f_1, \wh f_2)) ,\quad \wh f^3(\p A(\wh f_1, \wh f_2)) ,\quad\wh f^2( \p A(\wh f_1, \wh f_3)),\quad\wh f^3( \p A(\wh f_1, \wh f_3)),\\
& \wh f^4(\p A(\wh f_2, \wh f_3)),\quad  \wh f^5(\p A(\wh f_2, \wh f_3)),\quad  \wh f^6(\p A(\wh f_2, \wh f_3)),\quad  \wh f^7(\p A(\wh f_2, \wh f_3)).
\end{aligned}
}
\eeq
With the same  computations of \S \ref{sectA2}, we  compute all these coefficients and get that the equation  $   \langle \t, \p A \rangle = 0$ has the explicit expression
\beq
{\scalefont{0.90}
\begin{aligned} 
& \t^2_{12} \left( -  \frac{1}{\sqrt{6}} A^3_8  + \frac{1}{2\sqrt{3}} A^4_{10} + \frac{1}{2\sqrt{3}} A^6_{10}\right) +\t^3_{12}\left(\frac{1}{\sqrt{6}} A^2_8  + \frac{1}{2\sqrt{3}} A^5_{10} - \frac{1}{2\sqrt{3}} A^7_{10}\right)  +\\
 &+\t^2_{13} \left( -  \frac{1}{\sqrt{6}} A^3_9  + \frac{1}{2\sqrt{3}} A^5_{10} + \frac{1}{2\sqrt{3}} A^7_{10}\right) + \t^3_{13} \left(\frac{1}{\sqrt{6}} A^2_9  - \frac{1}{2\sqrt{3}} A^4_{10} + \frac{1}{2\sqrt{3}} A^6_{10}\right) +\\
& +\t^4_{23} \left( \frac{1}{2\sqrt{3}} A^2_9  +  \frac{1}{2\sqrt{3}} A^3_8 + \frac{1}{\sqrt{6}} A^4_{10}\right)  + \t^5_{23}\left(  \frac{1}{2\sqrt{3}} A^3_9  -  \frac{1}{2\sqrt{3}} A^2_8 + \frac{1}{\sqrt{6}} A^5_{10}\right) +\\
 & +\t^6_{23} \left( \frac{1}{2\sqrt{3}} A^2_9  -  \frac{1}{2\sqrt{3}} A^3_8 + \frac{1}{\sqrt{6}} A^6_{10}\right) +  \t^7_{23} \left( -  \frac{1}{2\sqrt{3}} A^3_9  -  \frac{1}{2\sqrt{3}} A^2_8 - \frac{1}{\sqrt{6}} A^7_{10}\right) = 0.
\end{aligned}
}
\eeq
Since this has to be satisfied for each $A$, factoring the components of $A$ we get that $\t \in  \ker \p^*|_{C^2_2(\gm_, \gg)}$ if and only if its components 
satisfy the system
\beq\nonumber
{\scalefont{0.90}
\begin{aligned}
 & \t^3_{12}   -   \frac{1}{\sqrt{2}} \t^5_{23}  -   \frac{1}{\sqrt{2}}   \t^7_{23}  = 0  ,\ \
  \t^2_{12}      -\frac{1}{\sqrt{2}} \t^4_{23}   +  \frac{1}{\sqrt{2}} \t^6_{23}    = 0, \\
 & \t^3_{13}   +\frac{1}{\sqrt{2}} \t^4_{23}  + \frac{1}{\sqrt{2}}\t^6_{23}  = 0,\ \ 
    \t^2_{13}   -  \frac{1}{\sqrt{2}}  \t^5_{23}   +   \frac{1}{\sqrt{2}}  \t^7_{23}   = 0,\\
  &  \t^4_{23}  + \frac{1}{\sqrt{2}}  \t^2_{12}   -\frac{1}{\sqrt{2}}   \t^3_{13}    = 0,\ \ 
    \t^5_{23} + \frac{1}{\sqrt{2}}  \t^3_{12}    + \frac{1}{\sqrt{2}} \t^2_{13}         = 0,\\
   &  \t^6_{23}  +\frac{1}{\sqrt{2}}   \t^2_{12}   +   \frac{1}{\sqrt{2}}  \t^3_{13}      = 0,\ \ 
     \t^7_{23}  + \frac{1}{\sqrt{2}}   \t^3_{12}  -  \frac{1}{\sqrt{2}} \t^2_{13} = 0.
\end{aligned} 
}
\eeq
A simple check shows that  this system has just the trivial solution. This   means that $\ker \p^*|_{C^2_2(\gm_, \gg)} = 0$ and that   (7.4) of \cite{MS} is equivalent  to 
$c^2_K  = 0$.
\par
\medskip
\subsection{The space $\ker \p^*|_{C^2_3(\gm_-, \gg)}$}
In this section we  determine explicitly the bilinear maps in 
$\Hom(\L^2 \gm_{-} , \gg)$ of shifting degree $+3$  that are in  $\ker \p^*$.  \par
A bilinear map  $\t \in \Hom(\L^2 \gm_-, \gg)$  of shifting degree $+3$ has  the form 
\beq
 \label{7.1}
 {\scalefont{0.90}
 \begin{aligned}
 &\t = \t^4_{12} f_4 \otimes (f^1 \wedge f^2) + \t^5_{12} f_5 \otimes (f^1 \wedge f^2) + \t^6_{12} f_6 \otimes (f^1 \wedge f^2) {+} \t^7_{12} f_7 \otimes (f^1 \wedge f^2)+\\
& + \t^4_{13} f_4 \otimes (f^1 \wedge f^3) + \t^5_{13} f_5 \otimes (f^1 \wedge f^3) + \t^6_{13} f_6 \otimes (f^1 \wedge f^3) + \t^7_{13} f_7 \otimes (f^1 \wedge f^3)+\\
 &+\t^8_{23} f_8 \otimes (f^2 \wedge f^3) +  \t^9_{23} f_9 \otimes (f^2 \wedge f^3) 
 .\\[-13pt]
 \end{aligned}
 }
 \eeq
 \ \\[-3pt]
As in the previous sections, this element is in $\ker \p^*$ if and only if  $\langle \p^* \t, A \rangle = -  \langle \t, \p A \rangle = 0$
for any  $A = A^\a_\b f_\a \otimes f^\b \in \Hom(\gh, \gg)$. From \eqref{7.1}, this condition is a linear equation in the 
components $\t^i_{jk}$   with coefficients
\beq\nonumber
{\scalefont{0.90}
\begin{aligned}
&\wh f^4(\p A(\wh f_1, \wh f_2)) ,\quad \wh f^5(\p A(\wh f_1, \wh f_2)) ,\quad \wh f^6(\p A(\wh f_1, \wh f_2)) ,\quad \wh f^7(\p A(\wh f_1, \wh f_2)), \\
&\wh f^4(\p A(\wh f_1, \wh f_3)) ,\quad \wh f^5(\p A(\wh f_1, \wh f_3)) ,\quad \wh f^6(\p A(\wh f_1, \wh f_3)) ,\quad \wh f^7(\p A(\wh f_1, \wh f_3)), \\
 &\wh f^8(\p A(\wh f_2, \wh f_3)),\quad  \wh f^9(\p A(\wh f_2, \wh f_3)).
 \end{aligned}
 }
 \eeq
We  explicitly compute these coefficients with the same standard computations  indicated  in  \S \ref{sectA2}. Then we obtain that the condition 
$\langle \t, \p A \rangle = 0$ has the explicit form  
\beq
{\scalefont{0.90}
\begin{aligned}
&\t^4_{12}\left(- \frac{1}{2 \sqrt{3}} A^2_{10}\right)+ \t^5_{12} \left(  - \frac{1}{2 \sqrt{3}} A^3_{10}\right) + \t^6_{12} \left(\frac{1}{\sqrt{3}} A^1_8 - \frac{1}{2\sqrt{3}} A^2_{10}\right)  +\\
&+ \t^7_{12} \left( \frac{1}{2\sqrt{3}} A^3_{10}\right) +  \t^4_{13}\left(\frac{1}{2 \sqrt{3}} A^3_{10}\right)  + \t^5_{13} \left(- \frac{1}{2 \sqrt{3}} A^2_{10}\right) +\\
 &+ \t^6_{13} \left(\frac{1}{\sqrt{3}} A^1_9 - \frac{1}{2 \sqrt{3}} A^3_{10}\right)  + \t^7_{13} \left(  - \frac{1}{2 \sqrt{3}} A^2_{10}\right) +\\
 &+\t^8_{23} \left(  \frac{1}{\sqrt{6}} A^1_8 - \frac{1}{\sqrt{6}} A^2_{10}\right)  +  \t^9_{23} \left(\frac{1}{\sqrt{6}} A^1_9  - \frac{1}{\sqrt{6}} A^3_{10}\right)  = 0.
\end{aligned}
}
\eeq
Since this needs to hold for each $A$, factoring the components of $A$ we get that $\t \in  \ker \p^*|_{C^2_3(\gm_, \gg)}$ if and only if its components 
satisfy the system for 
\beq\nonumber
{\scalefont{0.90}
\begin{aligned}
&   \frac{1}{2 \sqrt{3}} \t^4_{12}  + \frac{1}{2\sqrt{3}}  \t^6_{12}  + \frac{1}{2 \sqrt{3}} \t^5_{13}    + \frac{1}{2 \sqrt{3}}   \t^7_{13}    - \frac{1}{\sqrt{6}} \t^8_{23} = 0 ,  \\
 &    - \frac{1}{2 \sqrt{3}} \t^5_{12}  + \frac{1}{2\sqrt{3}}  \t^7_{12}  + \frac{1}{2 \sqrt{3}}  \t^4_{13}    - \frac{1}{2 \sqrt{3}} \t^6_{13}  + \frac{1}{\sqrt{6}} \t^9_{23}   = 0, \\
 & \frac{1}{\sqrt{3}}  \t^6_{12}   + \frac{1}{\sqrt{6}}  \t^8_{23}  = 0,\qquad
\frac{1}{\sqrt{3}}  \t^6_{13} + \frac{1}{\sqrt{6}} \t^9_{23}  = 0.
\end{aligned} 
}
\eeq
Using the last two equations to simplify the first two, the system reduces to 
\beq \label{eccoqua}
\begin{aligned}
&  \t^4_{12}  + \t^5_{13}       =   - 3\t^6_{12} -     \t^7_{13} , 
 \qquad   \t^4_{13}    -   \t^5_{12}       =   3 \t^6_{13} -   \t^7_{12} ,\\
 & \t^8_{23}  = - \sqrt{2}\,   \t^6_{12}  ,\qquad
    \t^9_{23}  =  - \sqrt{2}\, \t^6_{13} .
\end{aligned} 
\eeq
This means that  the space $\ker \p^*|_{C^2_2(\gm_, \gg)}$ is $6$-dimensional and that condition (8.4) of \cite{MS}  correspond to a system of linear equations on the curvature components
$$T^{0(10)}_{-2 \,-1(10)}\, \ \ T^{0(10)}_{-2 \,-1(01)},\ \ R^{0(10)}_{-2 \,-1(10)} ,\ \  R^{0(10)}_{-2 \,-1(01)},\ \ R^{1(10)}_{-1(10) \,-1(01)}.$$
In order to make explicit these equations, we have to convert the system \eqref{eccoqua} on the components of $\t$ in the basis $\cB$ into a system 
on the components of $\t$ in the standard  CR basis $\cB^{CR}$. For this purpose, we recall that 
\beq
\nonumber
{\scalefont{0.90}
\begin{aligned}
& e_{-2} =\sqrt{6} f_1,\quad
 e_{-1(10)}= \frac{\sqrt{6}}{2} (f_2 - i f_3),\quad e_{-1(01)}= \frac{\sqrt{6}}{2} (f_2 + i f_3),\\
& e_{0(10)}= \frac{\sqrt{12}}{2} (f_4 - i f_5),\quad e_{0(01)}= \frac{\sqrt{12}}{2} (f_4 + i f_5),\quad E_{0(10)}= \frac{\sqrt{12}}{2} (f_6 - i f_7),\\
&  e_{0(01)}= \frac{\sqrt{12}}{2} (f_6 + i f_7),\quad  E_{1(10)}= \frac{\sqrt{6}}{2} (f_8 - i f_9),\quad E_{1(01)}= \frac{\sqrt{6}}{2} (f_8 + i f_9)
\end{aligned}
}
\eeq
and that, for the dual vectors, 
\beq
\nonumber
{\scalefont{0.90}
\begin{aligned}
& e^{0(10)} = \frac{1}{\sqrt{12}} (f^4 + i f^5),\quad e^{0(01)} = \frac{1}{\sqrt{12}} (f^4 - i f^5),\quad 
E^{0(10)} = \frac{1}{\sqrt{12}} (f^6 + i f^7),\\
&  E^{0(01)} = \frac{1}{\sqrt{12}} (f^6 - i f^7), \quad  E^{1(10)}= \frac{1}{\sqrt{6}} (f^8 + i f^9),\quad E^{1(01)}=  \frac{1}{\sqrt{6}} (f^8 - i f^9).
\end{aligned}
}
\eeq
From this we get that 
\beq
{\scalefont{0.90}
\begin{aligned}
& T^{0(10)}_{-2\,-1(10)} = \frac{\sqrt{3}}{2} (f^4 + i f^5)(\t(f_1,  f_2) - i \t(f_1,f_3)) =\\[-4pt]
&\hskip 5cm =  \frac{\sqrt{3}}{2}(\t^4_{12} + \t^5_{13}) + i \frac{\sqrt{3}}{2} (- \t^4_{13} + \t^5_{12} ),
\end{aligned}
}
\eeq
\beq
{\scalefont{0.90}
\begin{aligned}
& T^{0(10)}_{-2 \,-1(01)}  = \frac{\sqrt{3}}{2} (f^4 + i f^5)(\t(f_1,  f_2) + i \t(f_1,f_3)) =\\[-4pt]
&\hskip 5cm =  \frac{\sqrt{3}}{2}(\t^4_{12} - \t^5_{13}) + i \frac{\sqrt{3}}{2} ( \t^4_{13} + \t^5_{12} ),
\end{aligned}
}
\eeq
\beq
{\scalefont{0.90}
\begin{aligned}
& R^{0(10)}_{-2\,-1(10)} = \frac{\sqrt{3}}{2} (f^6 + i f^7)(\t(f_1,  f_2) - i \t(f_1,f_3)) =\\[-4pt]
&\hskip 5cm =  \frac{\sqrt{3}}{2}(\t^6_{12} + \t^7_{13}) + i \frac{\sqrt{3}}{2} (- \t^6_{13} + \t^7_{12} ),
\end{aligned}
}
\eeq
\beq
{\scalefont{0.90}
\begin{aligned}
& R^{0(10)}_{-2 \,-1(01)}  = \frac{\sqrt{3}}{2} (f^6 + i f^7)(\t(f_1,  f_2) + i \t(f_1,f_3)) =\\[-4pt]
&\hskip 5cm =  \frac{\sqrt{3}}{2}(\t^6_{12} - \t^7_{13}) + i \frac{\sqrt{3}}{2} ( \t^6_{13} + \t^7_{12} ),
\end{aligned}
}
\eeq
 \beq
{\scalefont{0.90}
 R^{1(10)}_{-1(10) \,-1(01)}  =  \frac{\sqrt{3}}{\sqrt{2}} (f^8 + i f^9)( i \t(f_2 ,  f_3) ) =  \sqrt{ \frac{3}{2}}(  \t^8_{23} + i \t^9_{23}).
}
\eeq
From this we see that 
$$ {\scalefont{0.90} \frac{1}{\sqrt{3}} (\overline{R^{0(10)}_{-2\,-1(10)}}+ R^{0(10)}_{-2 \,-1(01)}) {=} \t^6_{12}+ i \t^6_{13},\  
 - i \frac{\sqrt{2}}{\sqrt{3}}R^{1(10)}_{-1(10) \,-1(01)} {=} \t^8_{23}+ i \t^9_{23},}$$
 which yields that  the last equation in \eqref{eccoqua} is equivalent to 
\beq\label{7.11}  {\scalefont{0.90} R^{1(10)}_{-1(10) \,-1(01)} = - i \overline{R^{0(10)}_{-2\,-1(10)}} - i R^{0(10)}_{-2 \,-1(01)}.} \eeq
On the other hand, since
\beq \nonumber
{\scalefont{0.90}
\begin{aligned}
& \frac{2}{\sqrt{3}}T^{0(10)}_{-2\,-1(10)}  = (\t^4_{12} + \t^5_{13}) - i  ( \t^4_{13} - \t^5_{12} ),\\[-1pt]
&\frac{2}{\sqrt{3}} \Bigg(\overline{R^{0(10)}_{-2\,-1(10)}} + R^{0(10)}_{-2\,-1(01)}\Bigg) +  \frac{2} {\sqrt{3}} \overline{R^{0(10)}_{-2 \,-1(10)}} =  3\t^6_{12} + \t^7_{13} + i (3\t^6_{13} -  \t^7_{12} ),
\end{aligned}
}
\eeq
we immediately see that the first two equations in \eqref{eccoqua} are equivalent to
\beq \label{7.12}  {\scalefont{0.90} T^{0(10)}_{-2\,-1(10)}  = - 2 \overline{ R^{0(10)}_{-2\,-1(10)}}  -   R^{0(10)}_{-2\,-1(01)}.}\eeq
Rearranging in an appropriate way the equations  \eqref{7.11} and \eqref{7.12},   we conclude  that condition (8.4) in \cite{MS} is equivalent to the following  equalities: 
\beq
{\scalefont{0.90}
\begin{aligned}
&    \overline{R^{0(10)}_{-2\,-1(10)} }   =   - \frac{1}{2} T^{0(10)}_{-2\,-1(10)}   - \frac{1}{2}   R^{0(10)}_{-2\,-1(01)}, \\[-4pt]
& R^{1(10)}_{-1(10) \,-1(01)} =  \frac{i}{2} T^{0(10)}_{-2\,-1(10)}  - \frac{i}{2} R^{0(10)}_{-2 \,-1(01)}.
\end{aligned}
}
\eeq
\par

\end{document}